\documentclass[11pt,a4paper]{amsart}

\usepackage[utf8]{inputenc}
\usepackage[T1]{fontenc}
\usepackage{xcolor}
\usepackage{listings}
\usepackage{amsmath,amssymb,amsthm,mathtools}
\usepackage{graphicx}

\usepackage{geometry}
\usepackage{placeins}
\usepackage{etoc}

\usepackage[hyphens]{url}
\usepackage[colorlinks=true, linkcolor=blue, citecolor=blue, urlcolor=blue]{hyperref}
\IfFileExists{orcidlink.sty}{\usepackage{orcidlink}}{\newcommand{\orcidlink}[1]{}}
\usepackage{doi}

\graphicspath{{./}{./figures/}}
\DeclareGraphicsExtensions{.pdf,.png,.jpg}

\hypersetup{
    pdftitle={Fej\'er--Kernel Prime Indicators},
    pdfauthor={Sebastian Fuchs (ORCID: 0009-0009-1237-4804)},
    pdfkeywords={Prime numbers, smooth prime indicator, analytic number theory, divisor function, sigma function, Fejér kernel, Fejér polynomial, truncation bounds, smooth lifts},
    pdfsubject={Primary 11A41; Secondary 11A25, 11N05, 11N99, 42A10}
}

\newcommand{\R}{\mathbb{R}}
\newcommand{\ii}{\mathrm{i}} 

\makeatletter
\renewcommand\paragraph{\@startsection{paragraph}{4}{\z@}%
  {1ex \@plus .2ex \@minus .2ex}%
  {0.8ex \@plus .2ex}%
  {\normalfont\bfseries}}
\makeatother

\theoremstyle{plain}
\newtheorem{theorem}{Theorem}[section]
\newtheorem{conjecture}[theorem]{Conjecture}
\newtheorem{lemma}[theorem]{Lemma}
\newtheorem{proposition}[theorem]{Proposition}
\newtheorem{corollary}[theorem]{Corollary}
\theoremstyle{definition}
\newtheorem{definition}[theorem]{Definition}
\theoremstyle{remark}
\newtheorem{remark}[theorem]{Remark}

\title[Fej\'er--Kernel Prime Indicators]{Fej\'er--Kernel Prime Indicators}
\author{Sebastian Fuchs \orcidlink{0009-0009-1237-4804}}
\address{Institut für Informatik, Humboldt-Universität zu Berlin, 10099 Berlin, Germany}
\email{\href{mailto:sebastian.fuchs@hu-berlin.de}{sebastian.fuchs@hu-berlin.de}}
\date{October 15, 2025}

\etocsetnexttocdepth{1}

\etocsettocstyle
  {\section*{\contentsname}\vspace{-0.4\baselineskip}\noindent\rule{\textwidth}{0.4pt}\par\medskip}
  {}

\etocsetstyle{section}
  {\par\etocskipfirstprefix}
  {\leftskip=1.4em \parindent=0pt \etocname{} \dotfill\etocpage\par}
  {}{}
\etocsetstyle{subsection}
  {}
  {\leftskip=3.0em \parindent=0pt \etocname{} \dotfill\etocpage\par}
  {}{}
\etocsetstyle{subsubsection}
  {}
  {\leftskip=5.0em \parindent=0pt \etocname{} \dotfill\etocpage\par}
  {}{}

\begin{document}

\maketitle

\begingroup\small
\noindent\textbf{Note:} A companion manuscript develops the divisor-filter mechanism into a general framework, the \href{https://arxiv.org/abs/2509.12297}{Fej\'er--Dirichlet Lift}, and establishes its connection to Dirichlet series and L-functions.
\par\medskip
\endgroup

\begin{abstract}
A $C^1$ prime indicator $\mathcal{P}\colon\mathbb{R}\to\mathbb{R}$ is constructed by applying the Fejér identity to the sine–quotient encoder of trial division. For integers $n\ge 2$, $\mathcal P(n)=0$ holds exactly for odd primes; $\mathcal P(2)>0$. For all non-integers $x>1$ one has $\mathcal P(x)>0$. The function is piecewise $C^\infty$ and its second derivative has jumps precisely at the squares $m^2$, with explicit sizes. Replacing the sharp cut-off by a smooth transition yields $C^\infty$ analogues $\mathcal{P}_\tau$ and $\mathcal{P}_\sigma$ with integer limits $\mathcal{P}_\tau(n;\kappa)\to \tau(n)-2$ and $\mathcal{P}_\sigma(n;\kappa)\to \sigma(n)-n-1$ as $\kappa\to\infty$, obtained from locally uniform convergence of derivative series. For large $\kappa$, numerical evidence indicates companion zeros near odd primes for $\mathcal{P}_\tau$ and an asymmetric pair for $\mathcal{P}_\sigma$. No assertion is made beyond integer input, and no statements are claimed about the prime number theorem or zero distributions of $L$-functions. The appendix includes two illustrative prime-counting sums.
\end{abstract}

\clearpage
\begingroup
\hypersetup{linkcolor=black} 
\tableofcontents
\endgroup
\hypersetup{linkcolor=blue}
\bigskip
\clearpage

\section{Introduction}

\subsection*{Background}
Approaches to the identification of the prime set may be grouped into four strands.
(i) \emph{Prime–representing formulae:} exact identities such as Mills/Willans/Wilson certify primality a posteriori but are computationally ineffective at scale.
(ii) \emph{Analytic zero–set constructions:} by Weierstrass factorization, entire functions with prescribed discrete zero sets exist, including ones vanishing exactly on primes; effectivity between integers is not provided~\cite{boas1954}. Related optical/scattering constructions appear in~\cite{petersen2019,li2019-optical,zhang2018-structure,torquato2018-scattering}.
(iii) \emph{Elementary trigonometric/kernel encoders:} divisibility can be encoded by quotients of sines
\[
Q(x,i):=\frac{\sin^2(\pi x)}{\sin^2(\pi x/i)}\qquad(i\in\mathbb N,\ i\ge2),
\]
which are continuous with removable singularities at integer arguments and relate to Fejér-type kernels via the Chebyshev identity \(\sin(ny)=U_{n-1}(\cos y)\sin y\)~\cite{katznelson2004,zygmund2002}.
(iv) \emph{L-function and sieve methods:} classical analytic number theory studies primes via \(\Lambda\), explicit formulas, zero-density estimates, and sieve weights; see~\cite{iwaniec2004,montgomery2007}.  

The construction in strand (iii) is adopted in the form of Fejér cosine polynomials to regularize trial division into a pointwise-defined function on $\mathbb{R}$. The resulting function retains a direct arithmetic meaning at integers and admits explicit control of smoothness and the locations and sizes of derivative jumps.

\subsection*{Method and Definition}
The Fejér identity converts the quotient-of-sines representation into a cosine polynomial
\[
F(x,i)=i+2\sum_{k=1}^{i-1}(i-k)\cos\!\left(\frac{2\pi k x}{i}\right)\qquad(i\ge2),
\]
which agrees with \(Q(x,i)\) away from removable singularities. At integer arguments, the normalized term $F(n,i)/i^2$ acts as an exact \emph{divisor filter}, evaluating to $\mathbf{1}_{i\mid n}$. The aggregate indicator is then defined as
\[
\mathcal P(x)=\frac{1}{x}\sum_{i=2}^{\lceil\sqrt{x}\rceil} F(x,i)\qquad(x>1).
\]
The construction produces a \(C^1\) function on \((0,\infty)\) with a direct arithmetic interpretation at integers:
\[
\mathcal P(n)=\frac{1}{n}\sum_{\substack{2\le d\le \lceil\sqrt n\rceil\\ d\mid n}} d^2.
\]
A resonant partial-fraction form and quantitative truncation bounds are provided in (Proposition~\ref{prop:RPF}, Theorem~\ref{thm:RPF-trunc}).

\subsection*{Main Results}
\begin{itemize}
  \item \textbf{Exact prime-zero characterization at odd integers.} For integers \(n\ge2\), \(\mathcal P(n)=0\) holds if and only if $n$ is an odd prime, while \(\mathcal P(2)>0\). For all non-integer \(x>1\), \(\mathcal P(x)>0\). (Theorem~\ref{thm:primezero}.)
  \item \textbf{A precise characterization of smoothness.} The function \(\mathcal P\) is \(C^1\) on \((0,\infty)\) and piecewise \(C^\infty\). Its second derivative has jump discontinuities exclusively at integer squares \(x=m^2\), governed by the explicit formula
  \[
  \Delta_{m^{2}}\mathcal{P}''=\frac{2\pi^{2}}{m^{2}\sin^{2}\!\bigl(\pi/(m+1)\bigr)}
  =2+\frac{4}{m}+O(m^{-2}).
  \]
  (Proposition~\ref{prop:second}.)
  \item \textbf{Extension to smooth analogues of arithmetic functions.} The construction is generalized to obtain \(C^\infty\) analogues \(\mathcal P_\tau(\cdot;\kappa)\) and \(\mathcal P_\sigma(\cdot;\kappa)\), which admit locally uniform control of all derivative series. For integers \(n\ge2\) and in the limit \(\kappa\to\infty\), these functions converge to values determined by classical arithmetic functions:
  \[
  \mathcal{P}_\tau(n;\kappa)\to \tau(n)-2,\qquad
  \mathcal{P}_\sigma(n;\kappa)\to \sigma(n)-n-1.
  \]
\end{itemize}

\subsection*{Scope and Limitations}
All statements about primes and composites concern integer input. No claims are made regarding the prime number theorem, densities, or zero-free regions. The indicator $\mathcal P$ characterizes odd primality at integers and requires $\Theta(\sqrt n)$ work, so no algorithmic advantage over trial division is implied. For non-integer input, only conjectural local features near odd primes are discussed for the smooth variants.

\subsection*{Computational Aspects}
Numerical experiments for the figures were produced by straightforward \(\Theta(\sqrt{x})\) evaluation using quotient-of-sines with resonance guards; the repository contains the scripts and parameter files required to reproduce the figures.\footnote{Repository: \url{https://github.com/SebastianFoxxx/analytic-prime-indicator}}


\section{Analytic construction of $\mathcal{P}$}

\paragraph*{Notation and disambiguation.}
The symbol $\ii=\mathrm{i}$ denotes the imaginary unit. Integer indices are consistently written as $i$; the context will always disambiguate from the imaginary unit. For lattices generated by an index, the explicit set notation $\{\, i k : k\in\mathbb{Z}\,\}$ is used instead of $i\mathbb{Z}$ to avoid ambiguity. Asymptotic symbols $O(\cdot)$ and $o(\cdot)$ have their standard meaning; unless stated otherwise, bounds are uniform on compact $x$–intervals $K\subset(0,\infty)$ for fixed integer parameters.

\subsection{Motivation}\label{subsec:motivation}
The construction uses an analytic representation of trial division. For a real number $x$ and an integer $i \ge 2$, divisibility can be expressed by the quotient
\[
Q(x,i) = \frac{\sin^2(\pi x)}{\sin^2(\pi x/i)}.
\]
The numerator vanishes if and only if $x$ is an integer. If $x$ is an integer, the denominator vanishes if and only if $i$ divides $x$. The quotient thus becomes indeterminate of the form $0/0$ precisely at integer arguments $x=n$ where $i$ is a divisor of $n$.

The indeterminacy is removed by the corresponding Fejér cosine polynomial $F(x,i)$. By Fejér’s identity the polynomial equals the quotient away from the removable singularities, which provides a well-defined extension at those points. For integers $n$, one has $F(n,i)=i^2$ if $i\mid n$ and $F(n,i)=0$ otherwise.

\subsection{Regularisation}
By the Fejér kernel identity~\cite{zygmund2002}
\[
\frac{\sin^{2}(r y)}{\sin^{2}y}=r+2\sum_{k=1}^{r-1}(r-k)\cos(2k y)\qquad (r\in\mathbb N),
\]
and with \(r=i\), \(y=\pi x/i\), one has the pointwise identity
\[
F(x,i)=i+2\sum_{k=1}^{i-1}(i-k)\cos\!\left(\frac{2\pi k x}{i}\right)
=\left(\frac{\sin(\pi x)}{\sin(\pi x/i)}\right)^{\!2},
\]
which is understood as the holomorphic continuation across the removable singularities. This yields a finite cosine polynomial that matches the quotient for \(x\notin\mathbb{Z}\) and reproduces the integer divisibility pattern.

\subsection{Definition}
The regularized function is constructed by summing the Fejér-kernel-based terms over all potential divisors up to $\sqrt{x}$.
Define $F(x,i)$ for $i\in\mathbb{N}$, $i\ge 2$, by
\begin{equation}\label{eq:Fxi}
F(x,i) = i+2\sum_{k=1}^{i-1}(i-k)\cos\left(\frac{2\pi k x}{i}\right).
\end{equation}


\begin{definition}\label{def:Px}
For $x>1$ the function $\mathcal{P}(x)$ is defined as
\begin{equation}\label{eq:Px}
\mathcal{P}(x)=\frac{1}{x}\sum_{i=2}^{\lceil\sqrt{x}\rceil} F(x,i),
\end{equation}
and $\mathcal{P}(x)=0$ is set for $x\le1$.
\end{definition}

\begin{remark}[On the choice of summation limit]
The ceiling function $\lceil\sqrt{x}\rceil$ is used instead of the more common floor function to ensure $C^1$-smoothness. As shown in Proposition~\ref{prop:C1}, at integer squares $x=m^2$, the new term $F(x, m+1)$ and its first derivative vanish, which guarantees a smooth transition. The floor function would introduce a discontinuity in the first derivative.
\end{remark}

Write $N(x)=\lceil\sqrt{x}\rceil$ for the upper summation limit; then \eqref{eq:Px} is $\mathcal{P}(x)=\frac{1}{x}\sum_{i=2}^{N(x)} F(x,i)$ for $x>1$.

\begin{remark}[Integer evaluation]\label{rem:integer-eval}
For an integer $n\ge 2$, Proposition~\ref{prop:fejer-properties} implies that $F(n,i)=i^2$ if $i\mid n$ and $F(n,i)=0$ otherwise. The definition in \eqref{eq:Px} therefore yields the exact evaluation
\[
\mathcal P(n)=\frac{1}{n}\sum_{\substack{d\mid n \\ 2\le d\le \lceil\sqrt{n}\rceil}} d^2.
\]
This sum vanishes exactly when $n$ has no divisors in the summation range. For $n>2$ this characterizes primality, but a direct evaluation requires $\Theta(\sqrt n)$ work and thus entails no algorithmic advantage over elementary trial division. At integer squares $n=m^2$, the limit $\lceil\sqrt{n}\rceil=m$ coincides with $\lfloor\sqrt{n}\rfloor$.
\end{remark}

\begin{proposition}[Properties of the Fejér term]\label{prop:fejer-properties}
For an integer $i \ge 2$, the function $F(\cdot,i)$ defined in \eqref{eq:Fxi} is an entire function with the following properties:
\begin{enumerate}
    \item[(a)] \textbf{Closed form and entire extension.} The trigonometric polynomial $F(\cdot,i)$ is entire. The map
    \[
    z\ \longmapsto\ \left(\frac{\sin(\pi z)}{\sin(\pi z/i)}\right)^{\!2}
    \]
    is meromorphic on $\mathbb{C}$ with removable singularities at the lattice points $z=ik$ ($k\in\mathbb{Z}$), because numerator and denominator vanish to the same order there. On the punctured plane $\mathbb{C}\setminus\{ik:k\in\mathbb{Z}\}$ the Fejér identity (with $y=\pi z/i$) yields pointwise equality with the cosine polynomial $F(\cdot,i)$. Removing the removable singularities produces the unique entire extension of this meromorphic function; by the identity theorem the resulting entire function coincides identically with $F(\cdot,i)$ on $\mathbb{C}$.
    \item[(b)] \textbf{Integer values.} For any $n \in \mathbb{Z}$,
    \[
    F(n,i) \;=\; \begin{cases} i^2, & \text{if } i \mid n, \\ 0, & \text{if } i \nmid n. \end{cases}
    \]
    \item[(c)] \textbf{Positivity on $\mathbb{R}$.} For all $x\in\mathbb{R}$, $F(x,i) \ge 0$, with strict inequality for $x\notin\mathbb{Z}$.
\end{enumerate}
\end{proposition}

\begin{proof}
The trigonometric polynomial in \eqref{eq:Fxi} defines an entire function. For $z \notin \{\,i k: k\in\mathbb{Z}\,\}$, Fejér's identity gives
\[
F(z,i)=\left(\frac{\sin(\pi z)}{\sin(\pi z/i)}\right)^{\!2}.
\]
The right-hand side is holomorphic on $\mathbb{C}\setminus\{\,i k: k\in\mathbb{Z}\,\}$ and has removable singularities at the points $z=i k$. Removing these poles yields an entire function that agrees with the trigonometric polynomial on a nonempty open set; the identity theorem then gives global equality.

For (b), if $i \nmid n$ then $\sin(\pi n/i)\ne 0$ while $\sin(\pi n)=0$, hence $F(n,i)=0$. If $i \mid n$, write $n=i m$ and apply L’Hôpital's rule to obtain $F(n,i)=i^2$.

For (c), if $x\notin\mathbb{Z}$ then $\sin(\pi x/i)\ne 0$ for every $i\ge 2$, so $F(x,i)=\bigl(\sin(\pi x)/\sin(\pi x/i)\bigr)^2>0$. At $x\in\{i k\}$ the quotient has a removable singularity and the extension yields the stated nonnegativity, with $F(x,i)=i^2>0$ when $x=i k$.
\end{proof}

\begin{proposition}[Resonant partial-fraction identity and convergence]\label{prop:RPF}
Let $i\in\mathbb{N}$, $i\ge 2$, and $z\in\mathbb{C}$. On the punctured plane $\mathbb{C}\setminus\{\,i k: k\in\mathbb{Z}\,\}$ one has
\[
\left(\frac{\sin(\pi z)}{\sin(\pi z/i)}\right)^{\!2}
=\frac{i^{2}}{\pi^{2}}\,\sin^{2}(\pi z)\,\sum_{k\in\mathbb{Z}}\frac{1}{(z-i k)^{2}},
\]
and the series on the right converges locally uniformly on compact subsets of this domain, so the displayed equality holds pointwise on $\mathbb{C}\setminus\{i k:k\in\mathbb{Z}\}$. Both sides therefore represent the same meromorphic function with removable singularities at $z=i k$. Removing these singularities produces the unique entire extension; by the identity theorem this entire extension coincides identically with the Fejér cosine polynomial $F(z,i)$.
\end{proposition}

\begin{proof}
The classical expansion $\pi^{2}\csc^{2}(\pi u)=\sum_{k\in\mathbb{Z}}(u-k)^{-2}$ converges locally uniformly on compact subsets of $\mathbb{C}\setminus\mathbb{Z}$ (Whittaker–Watson, Chapter VII, §2, cosecant expansion). With $u=z/i$,
\[
\csc^{2}\!\Bigl(\frac{\pi z}{i}\Bigr)=\frac{i^{2}}{\pi^{2}}\sum_{k\in\mathbb{Z}}\frac{1}{(z-i k)^{2}}
\]
on $\mathbb{C}\setminus\{\,i k: k\in\mathbb{Z}\,\}$. For any compact $K\subset\mathbb{C}\setminus\{\,i k: k\in\mathbb{Z}\,\}$ there exist constants $C_K,c_K>0$ such that $|z-i k|\ge c_K\,|k|$ for all $z\in K$ and $|k|$ sufficiently large. Hence
\[
\sum_{|k|>K_0}\frac{1}{|z-i k|^{2}}\ \le\ \frac{C_K}{c_K^{2}}\sum_{|k|>K_0}\frac{1}{k^{2}},
\]
and the Weierstrass M-test yields uniform convergence on $K$. Multiplication by $\sin^{2}(\pi z)$ preserves local uniform convergence on compacta disjoint from $\{i k\}$ because $\sin^{2}(\pi z)$ is holomorphic and bounded on each such compact set. This gives the identity
\[
\left(\frac{\sin(\pi z)}{\sin(\pi z/i)}\right)^{\!2}
=\frac{i^{2}}{\pi^{2}}\,\sin^{2}(\pi z)\,\sum_{k\in\mathbb{Z}}\frac{1}{(z-i k)^{2}}
\]
on $\mathbb{C}\setminus\{\,i k: k\in\mathbb{Z}\,\}$. At the points $z=i k$ both sides have removable singularities (the quotient has a $0/0$ form and the series side is a bounded holomorphic factor times a square-summable tail). Removing these singularities yields the unique entire extension of the common meromorphic function. Since both sides agree with the trigonometric polynomial $F(z,i)$ on a nonempty open subset of $\mathbb{C}$, the identity theorem implies global equality with $F(z,i)$.
\end{proof}

\begin{definition}[Nearest-lattice notation]\label{def:nearest-lattice}
For $x\in\mathbb{R}$ and an integer index $i\ge2$, set $m:=\mathrm{round}(x/i)$ with the fixed convention $\mathrm{round}(u):=\lfloor u+\tfrac12\rfloor$ (ties to the larger integer), and write $t:=x-im$. Then $|t|\le i/2$. In particular, $\sin(\pi x)=\pm \pi t+O(t^{3})$ and hence
\[
F(x,i)=\left(\frac{\sin(\pi x)}{\sin(\pi x/i)}\right)^{\!2}
=i^{2}+O\!\bigl(i^{2}(\pi t)^{2}\bigr)\qquad(t\to0),
\]
which recovers $F(im,i)=i^{2}$ by the removable-singularity extension.
\end{definition}

\begin{lemma}[Nearest-pole dominance bounds]\label{lem:nearest-pole-bounds}
Fix $i\ge2$ and $x\in\mathbb{R}$. Set $m:=\mathrm{round}(x/i)$ and $t:=x-i m$. The Fej\'er term can be written as $F(x,i) = \frac{i^2}{\pi^2}\sin^2(\pi x)\,\Sigma(x,i)$, where the pole sum $\Sigma(x,i) = \sum_{k\in\mathbb{Z}}(x-i k)^{-2}$ satisfies
\[
\frac{1}{t^2} \;\le\; \Sigma(x,i) \;\le\; \frac{1}{t^2} + \frac{\pi^2}{i^2}.
\]
\end{lemma}

\begin{proof}
The representation follows from Proposition~\ref{prop:RPF}. The sum $\Sigma(x,i)$ can be written as
\[
\Sigma(x,i) = \frac{1}{t^2}+\sum_{r=1}^{\infty}\left[\frac{1}{(x-i(m+r))^{2}}+\frac{1}{(x-i(m-r))^{2}}\right].
\]
The lower bound is obtained by discarding the non-negative tail sum. For the upper bound, $|t|\le i/2$ by definition, so
\[
|x-i(m\pm r)| = |t\mp ri| \ge ri - |t| \ge (r-\tfrac12)i\qquad(r\ge1).
\]
This gives the tail bound
\[
\sum_{r=1}^{\infty}\left[\dots\right] \le \frac{2}{i^2}\sum_{r=1}^{\infty}\frac{1}{(r-\tfrac12)^2} = \frac{2}{i^2}\cdot\frac{\pi^2}{2} = \frac{\pi^2}{i^2},
\]
where the sum identity follows from the standard cosecant series evaluated at $z=1/2$.
\end{proof}

\begin{remark}[Integer arguments and removable singularities]
If $x=n\in\mathbb{Z}$, the inequalities in Lemma~\ref{lem:nearest-pole-bounds} are to be understood in the removable-singularity sense. If $i\nmid n$ then $\sin(\pi n)=0$ and $F(n,i)=0$, so both sides read $0\le 0$. If $i\mid n$, write $t=x-n$; since $\sin(\pi x)\sim \pi(x-n)$ as $x\to n$, one has
\[
\frac{i^{2}}{\pi^{2}}\frac{\sin^{2}(\pi x)}{t^{2}}\ \longrightarrow\ i^{2}\ =\ F(n,i),
\]
which shows that both bounds converge to the common value at $x=n$.
\end{remark}

\begin{remark}[An identity on the divisor lattice]
For integer arguments $n$, the normalized Fejér terms satisfy a multiplicative property on the divisor lattice. Since $F(n,r)/r^2$ is the indicator function $\mathbf{1}_{r\mid n}$ for integers, the identity $\mathbf{1}_{i\mid n}\mathbf{1}_{j\mid n}=\mathbf{1}_{\operatorname{lcm}(i,j)\mid n}$ implies
\[
\frac{F(n,i)}{i^2}\cdot \frac{F(n,j)}{j^2}\;=\;\frac{F\bigl(n,\operatorname{lcm}(i,j)\bigr)}{\operatorname{lcm}(i,j)^2}
\]
for any integers $i,j\ge 2$. This property is specific to integer inputs.
\end{remark}

\begin{remark}[Classical references]
Standard facts on Fejér kernels, cosine–polynomial representations, and Cesàro means are classical; see Zygmund~\cite[Chap.~I, §5; Chap.~III]{zygmund2002}, Katznelson~\cite[Chap.~I]{katznelson2004}, and Stein–Shakarchi~\cite[Chap.~2]{stein2003}. Nonnegativity follows from the Fejér–Riesz factorization for nonnegative trigonometric polynomials; cf.\ \cite[Chap.~V]{zygmund2002} and \cite[Chap.~II]{katznelson2004}.
\end{remark}


\section{Smoothness properties of $\mathcal{P}$}

\subsection{Global $C^{1}$-smoothness}

\begin{proposition}[$C^1$-Smoothness of $\mathcal{P}$]\label{prop:C1}
The function $\mathcal{P}$ is continuous on $\R$ and of class $C^1$ on $(0,\infty)$.
\end{proposition}

\begin{proof}
For $x\le 1$ one has $\mathcal P(x)\equiv 0$. On each open interval $((m-1)^2,m^2)$ the index $N(x)$ is constant, hence $\mathcal P$ is a finite sum of $C^\infty$ functions there. It remains to check the junctions $x=m^2$ and $x=1$.

\emph{Claim.} For every $n\in\mathbb Z$ and every $i\ge 2$, $F'(n,i)=0$. Indeed, write $F(x,i)=s(x)^2$ with $s(x):=\sin(\pi x)/\sin(\pi x/i)$. If $i\nmid n$ then $s(n)=0$ and hence $F'(n,i)=2s(n)s'(n)=0$. If $i\mid n$, write $n=im$ and $x=n+t$. Then
\[
s(n+t)=\frac{\sin(\pi(n+t))}{\sin(\pi(n+t)/i)}
=\frac{(-1)^n\sin(\pi t)}{(-1)^m\sin(\pi t/i)}
=(-1)^{n-m}\,\frac{\text{odd}(t)}{\text{odd}(t)}.
\]
The quotient of two odd functions is even in $t$, so $s'(n)=0$ and again $F'(n,i)=2s(n)s'(n)=0$.

At $x=m^2$ the only term changing is $G_m(x):=F(x,m{+}1)/x$. Since $(m{+}1)\nmid m^2$, $F(m^2,m{+}1)=0$ and, by the claim, $F'(m^2,m{+}1)=0$. For any $H$ with $H(x_0)=H'(x_0)=0$ and $x_0>0$,
\[
\Bigl(\frac{H}{x}\Bigr)'(x_0)=\frac{x_0H'(x_0)-H(x_0)}{x_0^2}=0.
\]
Thus $G_m'(m^2)=0$, and the one-sided first derivatives match. The same argument at $x=1$ with $G_1(x):=F(x,2)/x$ gives $G_1'(1)=0$. Therefore $\mathcal P$ is $C^1$ on $(0,\infty)$.
\end{proof}

\begin{lemma}[Heaviside junction lemma]\label{lem:heaviside-junction}
Let $x_0>0$. Suppose $H$ is $C^1$ on a neighborhood of $x_0$ and $C^2$ on each side of $x_0$, and $G\in C^2$ on a neighborhood of $x_0$ with $G(x_0)=G'(x_0)=0$. Define
\[
\widetilde H(x)=H(x)+G(x)\,\mathbf 1_{[x_0,\infty)}(x).
\]
Then $\widetilde H$ is $C^1$ at $x_0$, and
\[
\widetilde H''(x_0{+})-\widetilde H''(x_0{-})=G''(x_0).
\]

\begin{proof}
Since $G(x_0)=G'(x_0)=0$, both $\widetilde H$ and $\widetilde H'$ coincide from the left and from the right at $x_0$, hence $\widetilde H\in C^1$. On $(-\infty,x_0)$ one has $\widetilde H''=H''$, and on $(x_0,\infty)$ one has $\widetilde H''=H''+G''$. Therefore
\[
\widetilde H''(x_0{+})-\widetilde H''(x_0{-})=\bigl(H''(x_0)+G''(x_0)\bigr)-H''(x_0)=G''(x_0).
\]
\end{proof}
\end{lemma}

\subsection{Second-derivative jumps}

\begin{proposition}[Second-derivative jump at squares]\label{prop:second}
For $m\ge1$, the second derivative $\mathcal P''$ admits one-sided limits at $x=m^2$, and the jump
\[
\Delta_{m^2}\mathcal P''\ :=\ \mathcal P''(m^2{+})-\mathcal P''(m^2{-})
\]
is strictly positive and equals
\[
\Delta_{m^2}\mathcal{P}''=\frac{2\pi^2}{m^2\,\sin^2\!\bigl(\pi/(m+1)\bigr)}.
\]
In particular, for $m=1$ one obtains $\Delta_{1^2}\mathcal P''=2\pi^2$.
\end{proposition}

\begin{proof}
Across $x=m^2$ only the new index $i=m{+}1$ enters, so $G_m(x):=F(x,m{+}1)/x$ is the only additional term; the remaining finite sum
\[
H(x)=\frac{1}{x}\sum_{i=2}^{m}F(x,i)
\]
is $C^\infty$ on each open side and, since $1/x$ is smooth at $x=m^2>0$, also $C^1$ across $x=m^2$. Both $H$ and $G(x)=F(x,m{+}1)/x$ satisfy the hypotheses of Lemma~\ref{lem:heaviside-junction}: $H$ is $C^2$ on each side and $C^1$ in a neighborhood of $m^2$, while $G\in C^2$ near $m^2$ with $G(m^2)=G'(m^2)=0$ (cf.\ Proposition~\ref{prop:C1}). Hence
\[
\Delta_{m^2}\mathcal P''=G_m''(m^2).
\]
Since $\bigl(\tfrac{F}{x}\bigr)'=\frac{xF'-F}{x^2}$ and $\bigl(\tfrac{F}{x}\bigr)''=\frac{F''}{x}-\frac{2(xF'-F)}{x^3}$, the identities $F(m^2,m{+}1)=F'(m^2,m{+}1)=0$ imply
\[
G_m''(m^2)=\frac{F''(m^2,m{+}1)}{m^2}.
\]
Let $s(x)=\sin(\pi x)/\sin(\pi x/(m{+}1))$. Then $F=s^2$ and $F''(m^2,m{+}1)=2\,[s'(m^2)]^2$ because $s(m^2)=0$. Differentiation yields
\[
s'(x)=\frac{\pi\cos(\pi x)\sin(\pi x/(m{+}1))-\frac{\pi}{m{+}1}\sin(\pi x)\cos(\pi x/(m{+}1))}{\sin^2(\pi x/(m{+}1))}.
\]
At $x=m^2$ one has $\sin(\pi m^2)=0$ and, using $m^2=(m-1)(m{+}1)+1$,
\[
s'(m^2)=\frac{\pi\,(-1)^{m^2}}{\sin\!\bigl(\pi m^2/(m{+}1)\bigr)},
\qquad
\bigl|s'(m^2)\bigr|=\frac{\pi}{\sin\!\bigl(\pi/(m{+}1)\bigr)}.
\]
Since $F''(m^2,m{+}1)=2\,[s'(m^2)]^2$, the sign is immaterial.
Therefore $F''(m^2,m{+}1)=2\pi^2\,\csc^2\!\bigl(\pi/(m{+}1)\bigr)$ and
\[
\Delta_{m^2}\mathcal P''=\frac{2\pi^2}{m^2\,\sin^2\!\bigl(\pi/(m{+}1)\bigr)}.
\]
\end{proof}

\begin{remark}[Asymptotics of the jump size]
From
\[
\Delta_{m^2}\mathcal{P}''=\frac{2\pi^2}{m^2\,\sin^2\!\bigl(\pi/(m+1)\bigr)}
\]
set $z=\pi/(m+1)$ and use $\csc^2 z=z^{-2}+\frac{1}{3}+O(z^2)$ as $z\to 0$. Then
\[
\frac{1}{\sin^2\!\bigl(\pi/(m+1)\bigr)}
=\frac{1}{z^2}\Bigl(1+\tfrac{1}{3}z^2+O(z^4)\Bigr)
=(m+1)^2\Bigl(\tfrac{1}{\pi^2}+\tfrac{1}{3(m+1)^2}+O(m^{-4})\Bigr).
\]
With $(m+1)^2=m^2(1+\tfrac{2}{m}+\tfrac{1}{m^2})$ this gives
\[
\Delta_{m^2}\mathcal{P}''=2\Bigl(1+\tfrac{2}{m}+\tfrac{1}{m^2}\Bigr)\Bigl(1+\tfrac{\pi^2}{3(m+1)^2}+O(m^{-4})\Bigr)
=2+\frac{4}{m}+\frac{2+2\pi^2/3}{m^2}+O(m^{-3}).
\]
In particular, $\Delta_{m^2}\mathcal P''\downarrow 2$ as $m\to\infty$.
\end{remark}

\begin{figure}[!htbp]
\centering
\includegraphics[width=0.9\textwidth]{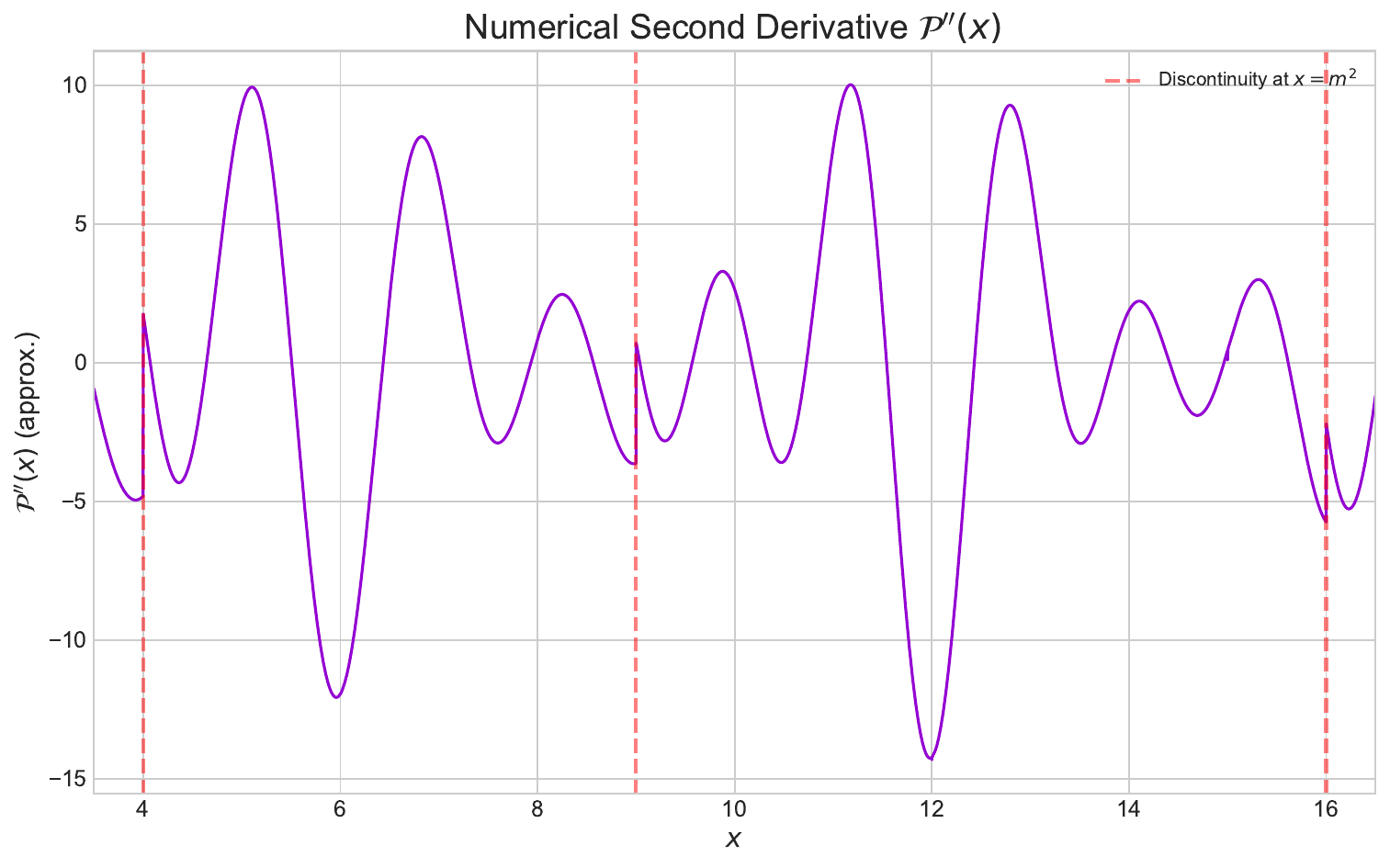}
\caption{Central finite-difference profile of $\mathcal{P}''(x)$ on $[3.5,16.5]$. Discontinuities at $x=4,9,16$ (integer squares) match Proposition~\ref{prop:second}; the observed jump heights agree with $\Delta_{m^2}\mathcal P''=2\pi^2\,[m^2\sin^2(\pi/(m+1))]^{-1}$.}
\label{fig:secondderivative}
\end{figure}

\FloatBarrier

\section{The Prime–Zero Property of $\mathcal{P}(x)$}

\begin{remark}
By construction, $\mathcal P(2)=2$. Hence the zero statement in Theorem~\ref{thm:primezero} concerns odd primes only.
\end{remark}

\noindent The divisor-filter identity $F(n,i)/i^2=\mathbf{1}_{i\mid n}$ is stated in Proposition~\ref{prop:fejer-properties}\,(b) and is used implicitly below.

\begin{theorem}[Odd prime zeros and positivity away from integers]\label{thm:primezero}
For all real $x>1$, the function $\mathcal P(x)$ is strictly positive whenever $x\notin\mathbb Z$. Moreover,
\[
\mathcal P(x)=0\quad\Longleftrightarrow\quad x\in\mathbb N\ \text{and $x$ is an odd prime}.
\]
In particular, $\mathcal P(2)=2$.
\end{theorem}
\begin{proof}
Let $x>1$ with $x\notin\mathbb Z$. Then $x/i\notin\mathbb Z$ for every integer $i\ge 2$; otherwise $x=i\cdot (x/i)\in\mathbb Z$, a contradiction. Hence, by Proposition~\ref{prop:fejer-properties}(c), $F(x,i)>0$ holds for all $i\ge 2$, so each summand in the finite sum is strictly positive and therefore $\mathcal P(x)>0$.

Let $x=n\in\mathbb N$, $n\ge 2$. If $n$ is an odd prime, then $i\nmid n$ for all $2\le i\le \lceil\sqrt n\rceil$, implying $F(n,i)=0$ and $\mathcal P(n)=0$. If $n$ is composite, there exists $d$ with $2\le d\le \sqrt n$ and $d\mid n$, hence $F(n,d)=d^2>0$ and $\mathcal P(n)>0$. Finally, at $n=2$ the sum reduces to the single index $i=2$ with $F(2,2)=4$, so $\mathcal P(2)=\frac{4}{2}=2>0$.
\end{proof}

\begin{lemma}[Local quadratic behavior near an odd prime]\label{lem:local-quadratic}
Let $p\ge 3$ be prime and put $N(p):=\lceil \sqrt{p}\rceil$. Define
\[
c_p\ :=\ \tfrac12\min_{\,2\le i\le N(p)}\bigl|\sin(\pi p/i)\bigr|\ >\ 0,
\qquad
\delta_p^{(1)}\ :=\ \frac{c_p}{\pi},
\qquad
\delta_p^{(2)}\ :=\ \tfrac12\min\bigl\{\,p-(N(p)-1)^2,\ N(p)^2-p\,\bigr\}.
\]
Let $\delta_p:=\min\{\delta_p^{(1)},\delta_p^{(2)}\}$. Then for all $x$ with $|x-p|<\delta_p$ one has $N(x)=N(p)$ and
\[
\mathcal P(x)=C_p\,(x-p)^2+O\bigl((x-p)^3\bigr),
\qquad
C_p=\frac{\pi^2}{p}\sum_{i=2}^{N(p)}\frac{1}{\sin^2(\pi p/i)}\ \ge\ \frac{\pi^2}{p}\,.
\]
The implicit constant in the $O(\cdot)$ term depends only on $p$ and the bound is uniform over all $i\le N(p)$.

\begin{proof}
For every $i\in\{2,\dots,N(p)\}$ one has $\sin(\pi p/i)\neq 0$ because $p$ is prime and $i<p$, hence $c_p>0$. Moreover,
\[
\Bigl|\partial_x \sin\!\Bigl(\frac{\pi x}{i}\Bigr)\Bigr|\le \frac{\pi}{2}\qquad(i\ge 2).
\]
Thus, for $|x-p|<\delta_p^{(1)}:=c_p/\pi$,
\[
\bigl|\sin(\tfrac{\pi x}{i})\bigr|\ge \bigl|\sin(\tfrac{\pi p}{i})\bigr|-\frac{\pi}{2}|x-p|\ge \frac{c_p}{2},
\]
uniformly in $2\le i\le N(p)$. Choosing $\delta_p^{(2)}$ as in the statement keeps $N(x)$ constant. With $\delta_p:=\min\{\delta_p^{(1)},\delta_p^{(2)}\}$, Taylor expansions yield
\[
\sin(\pi x)=(-1)^p\Bigl(\pi(x-p)-\frac{\pi^3}{6}(x-p)^3\Bigr)+O\bigl(|x-p|^5\bigr),
\]
\[
\sin\!\Bigl(\frac{\pi x}{i}\Bigr)=\sin\!\Bigl(\frac{\pi p}{i}\Bigr)+\frac{\pi}{i}\cos\!\Bigl(\frac{\pi p}{i}\Bigr)(x-p)+O\bigl((x-p)^2\bigr),
\]
where the $O(\cdot)$ is uniform in $i$ for $2\le i\le N(p)$. Division is legitimate since $|\sin(\pi x/i)|\ge c_p/2$ on $|x-p|<\delta_p$, and hence
\[
F(x,i)=\left(\frac{\sin(\pi x)}{\sin(\pi x/i)}\right)^{\!2}
=\frac{\pi^2}{\sin^2(\pi p/i)}(x-p)^2+R_i(x),
\]
with a uniform bound $|R_i(x)|\le C_p^{(0)}\,|x-p|^3$ for some constant $C_p^{(0)}$ depending only on $p$. Summation over $i\le N(p)$ preserves the $O(|x-p|^3)$ scale. Since $N(x)=N(p)$ and $1/x=1/p+O(x-p)$ on $|x-p|<\delta_p$,
\[
\mathcal P(x)=\frac{1}{x}\sum_{i=2}^{N(p)}F(x,i)
=\frac{1}{p}\sum_{i=2}^{N(p)}\frac{\pi^2}{\sin^2(\pi p/i)}(x-p)^2+O\bigl(|x-p|^3\bigr),
\]
with an explicit constant $C_p:=C_p^{(0)}+\tfrac{\pi^2}{p}\sum_{i\le N(p)}\bigl|\sin(\pi p/i)\bigr|^{-2}$ controlling the $O(\cdot)$ uniformly. The lower bound $C_p\ge \pi^2/p$ follows from $i=2$ and $\sin(\pi p/2)=\pm 1$ for odd $p$.
\end{proof}
\end{lemma}

\begin{remark}[Extension to all primes]
The indicator $\mathcal{P}(x)$ does not vanish at $p=2$. The companion manuscript~\cite{fuchs2025-fejer} develops entire-function analogues based on infinite weighted sums of divisor filters, which extend the prime-zero property to all primes, including $p=2$, at the cost of requiring a limiting process ($\kappa\to\infty$) for exact integer evaluation.
\end{remark}

\begin{figure}[!htbp]
\centering
\includegraphics[width=0.9\textwidth]{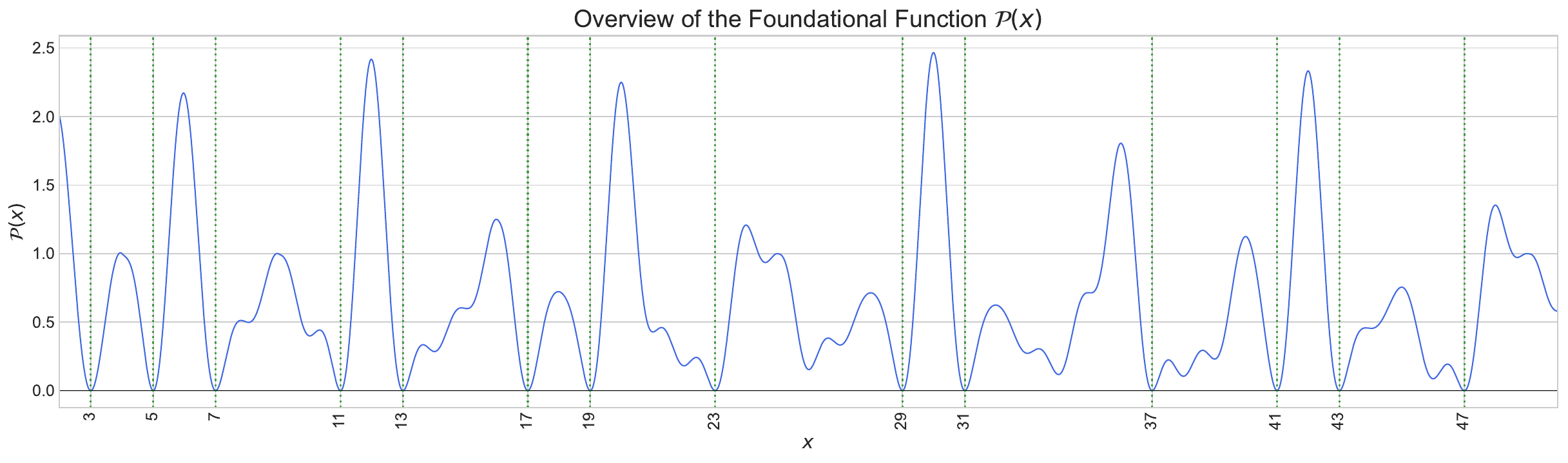}
\caption{Numerical profile of $\mathcal{P}(x)$ for $2\le x\le 50$.
Dotted vertical lines mark odd primes;
$\mathcal{P}(x)$ vanishes at these positions in accordance with Theorem~\ref{thm:primezero}.}
\label{fig:overview}
\end{figure}

\begin{figure}[!htbp]
\centering
\includegraphics[width=0.8\textwidth]{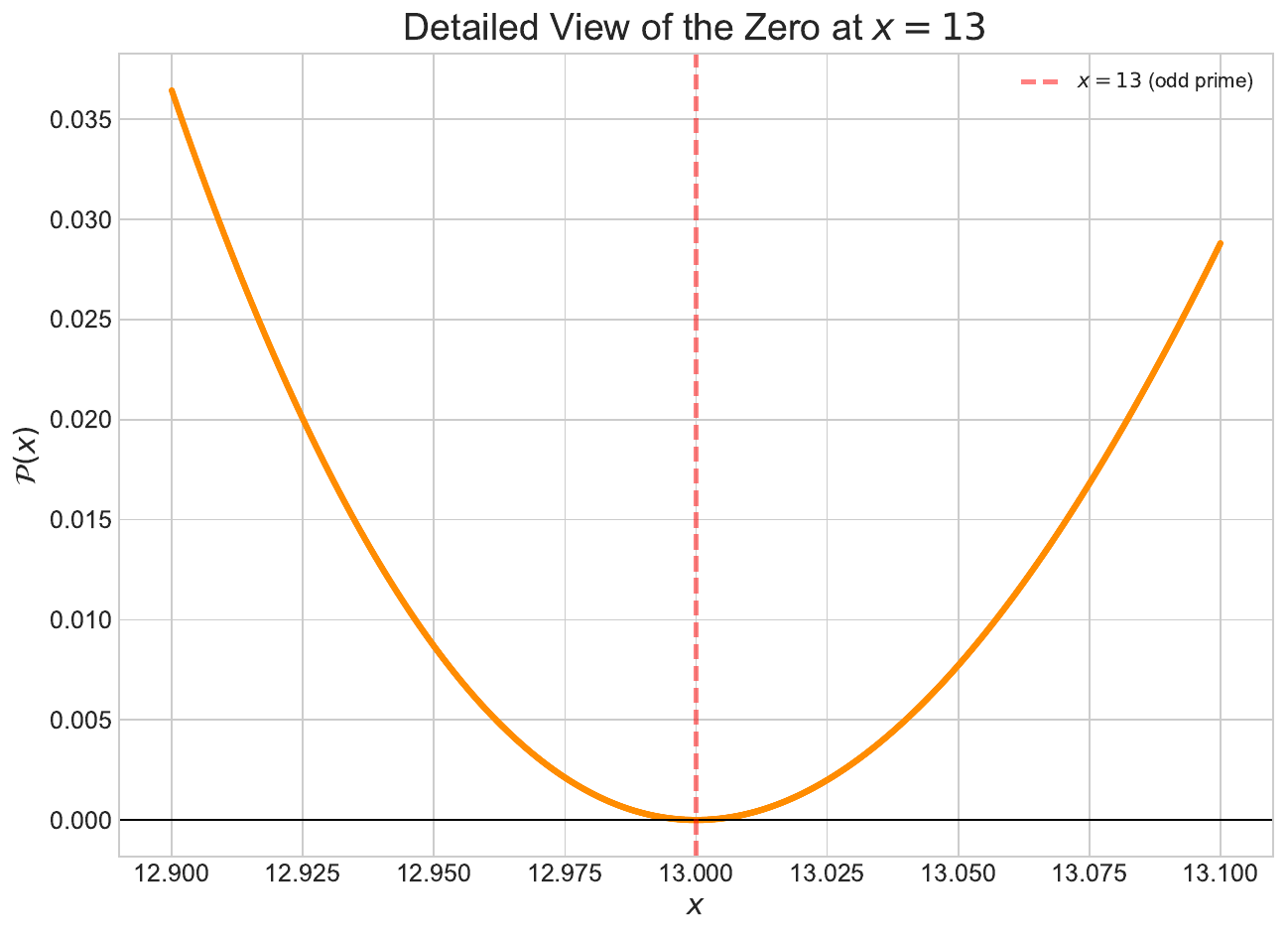}
\caption{Local behavior of $\mathcal{P}(x)$ near $x=13$ (an odd prime).
A zero occurs at $x=13$ without a sign change, consistent with nonnegativity and the odd prime–zero property.}
\label{fig:zoom13}
\end{figure}

\FloatBarrier

\section{Analytic Structure and Computational Complexity}

This section investigates deeper analytic properties of $\mathcal{P}(x)$, including a resonant partial-fraction identity that provides explicit error bounds for numerical evaluation schemes.

\subsection{Lower and upper bounds for composite integers}
A quantitative lower bound for composite $n$ follows from Theorem~\ref{thm:primezero}:
\[
\mathcal{P}(n) = \frac{1}{n} \sum_{\substack{d\mid n \\ 2 \le d \le \sqrt{n}}} d^2 \;\ge\; \frac{4}{n}.
\]
Equality holds precisely when the only divisor $d$ in the range $2\le d\le \sqrt{n}$ is $d=2$. This occurs if and only if $n=2p$ with odd prime $p$, or $n\in\{4,8\}$. For higher powers of two, $n=2^k$ with $k\ge 4$, the divisor $d=4$ satisfies $2\le d\le\sqrt{n}$ and contributes an additional positive term, and for $n=2p^a$ with $a\ge 2$ an extra divisor $d=p$ lies in the range; hence equality fails in all remaining composite cases. A complementary upper bound follows from
\[
\sum_{d\le \sqrt{n}} d^2 \;=\; \frac{1}{6}\bigl\lfloor\sqrt{n}\bigr\rfloor\bigl(\bigl\lfloor\sqrt{n}\bigr\rfloor+1\bigr)\bigl(2\bigl\lfloor\sqrt{n}\bigr\rfloor+1\bigr)
\;=\; \frac{1}{3}n^{3/2} + O(n),
\]
which implies
\[
\mathcal{P}(n) \;\le\; \frac{1}{n}\sum_{d\le \sqrt{n}} d^2
\;=\; \frac{1}{3}\,n^{1/2} + O(1).
\]
Thus $\mathcal{P}(n)$ lies between $4/n$ and $\ll n^{1/2}$, which determines the scale of $\mathcal{P}(n)$ on composite integers.

\begin{remark}[Equality case for the lower bound]
Equality $\mathcal{P}(n)=4/n$ holds iff the only divisor of $n$ in $[2,N(n)]$ is $2$; e.g. $n=2p$ with odd $p$, or $n\in\{4,8\}$.
\end{remark}

\subsection{Computational Complexity}

\smallskip
\noindent\textit{(A) Direct Fejér-sum evaluation.}
Using \eqref{eq:Fxi} inside \eqref{eq:Px} leads to $\sum_{i\le \lceil\sqrt{x}\rceil}\sum_{k<i}$ cosine evaluations.
This costs $\sum_{i\le N(x)}(i-1)=\tfrac12N(x)(N(x)-1)=\Theta(x)$ operations, thus $O(x)$.

\smallskip
\noindent\textit{(B) Quotient-of-sines evaluation.}
For $x\notin\mathbb{Z}$, each summand can be computed as $F(x,i)=\sin^2(\pi x)/\sin^2(\pi x/i)$, avoiding the inner sum.
This yields $\Theta(N(x))=\Theta(\sqrt{x})$ evaluations overall, thus $O(\sqrt{x})$.
At integer $x=n$ and divisors $i\mid n$, the limit is evaluated via L’Hôpital’s rule, giving $F(n,i)=i^2$; numerically this can be handled either by a branch for $|x/i-\mathrm{round}(x/i)|<\varepsilon$ or by evaluating $\lim_{t\to0}\sin^2(\pi n)/\sin^2(\pi(n/i+t))=i^2$.
For non-divisors, the denominator stays bounded away from $0$ and the quotient is safe.
The finitely many divisor cases do not alter the overall $O(\sqrt{x})$ complexity.

\smallskip
\noindent\textit{(C) Resonant partial fractions (RPF).}
By Proposition~\ref{prop:RPF},
\[
F(x,i)=\frac{i^{2}}{\pi^{2}}\sin^{2}(\pi x)\sum_{k\in\mathbb{Z}}\frac{1}{(x-ik)^{2}}.
\]
For evaluation, sum only the $(2K{+}1)$ nearest poles around $m=\mathrm{round}(x/i)$ (any fixed tie-breaking rule is admissible); Theorem~\ref{thm:RPF-trunc} yields an \emph{absolute} error bound $O(1/K)$, independent of $i$ apart from the prefactor $\sin^{2}(\pi x)$.  
Away from integers, Corollary~\ref{cor:RPF-relative} gives a \emph{relative} bound of size $\le \frac{1}{2(K+1/2)}$. With fixed small $K$ (e.g. $K=1,2$) the total cost remains $O(\sqrt{x})$, and stability near resonances is preserved without branch switching.

\paragraph*{Stable evaluation near resonances}
When $x/i$ is close to an integer, the quotient form loses significance. In the regime $\bigl|x/i-\operatorname{round}(x/i)\bigr|<\varepsilon$ use a constant-time local Taylor surrogate for 
\[
F(x,i)=\left(\frac{\sin(\pi x)}{\sin(\pi x/i)}\right)^{\!2}.
\]
Write $x=im+t$ with $m=\operatorname{round}(x/i)$ and $|t|<\varepsilon i$. Then
\[
F(x,i)=i^2\Bigl(1-\alpha_i(\pi t)^2\Bigr)^2+R_i(t),\qquad 
\alpha_i=\tfrac{1}{6}\Bigl(1-\tfrac{1}{i^2}\Bigr),
\]
with the explicit remainder bound from Lemma~\ref{lem:resonance}. At exact divisors ($t=0$) return $i^2$. This branch costs $O(1)$ per index, hence the total remains $O(\sqrt{x})$.

\medskip
\begin{lemma}[Resonant indices and local even expansion]\label{lem:resonance}
Let $\varepsilon\in(0,\tfrac14]$, $i\ge2$, and write $m=\mathrm{round}(x/i)$ and $x=mi+t$ with $|t|\le \varepsilon i$. Then
\[
\frac{\sin(\pi x)}{\sin(\pi x/i)} \;=\; (-1)^{m(i-1)}\,i\Bigl(1-\alpha_i(\pi t)^2\Bigr)\;+\;R_i^{(s)}(t),
\qquad
\alpha_i=\frac{1}{6}\Bigl(1-\frac{1}{i^{2}}\Bigr),
\]
and the remainder satisfies the uniform bound
\[
\bigl|R_i^{(s)}(t)\bigr|\;\le\; i\,C(\varepsilon)\,(\pi|t|)^4,
\qquad
C(\varepsilon):=\frac{1}{120}\,c(\varepsilon)^{-2},\quad
c(\varepsilon):=1-\frac{\pi^2\varepsilon^2}{6}-\frac{\pi^4\varepsilon^4}{120}.
\]
Consequently,
\[
F(x,i)=\left(\frac{\sin(\pi x)}{\sin(\pi x/i)}\right)^{\!2}
= i^2\Bigl(1-2\alpha_i(\pi t)^2\Bigr)\;+\;R_i^{(F)}(t),
\qquad
\bigl|R_i^{(F)}(t)\bigr|\;\le\; i^2\,\widetilde C(\varepsilon)\,(\pi|t|)^4
\]
for some explicit $\widetilde C(\varepsilon)\asymp C(\varepsilon)$.
Moreover, for $i\le\sqrt{x}$ and $m=\mathrm{round}(x/i)$, the resonance set $\{\,i:|x/i-m|<\varepsilon\,\}$ equals 
\[
I_m:=\Bigl(\frac{x}{m+\varepsilon},\,\frac{x}{m-\varepsilon}\Bigr)\cap (2,\infty),
\]
and its total length satisfies, for all $x\ge 4$,
\[
\sum_{m\ge \lceil \sqrt{x}\rceil}\! |I_m| 
\;=\;\sum_{m\ge \lceil \sqrt{x}\rceil}\frac{2\varepsilon\,x}{m^2-\varepsilon^2}
\;\le\; \frac{2\varepsilon\,x}{\lceil \sqrt{x}\rceil-1-\varepsilon}
\;\le\; 4\,\varepsilon\,\sqrt{x}.
\]
For $1\le x<4$ the equality remains valid; a finite bound suffices in applications and is not required explicitly below.

\begin{proof}
Taylor expansions with Lagrange remainder give, for $|t|\le\varepsilon i$,
\[
\sin(\pi t)=\pi t-\frac{(\pi t)^3}{6}+R_5(\pi t),\qquad |R_5(z)|\le \frac{|z|^5}{120},
\]
\[
\sin\!\Bigl(\frac{\pi t}{i}\Bigr)=\frac{\pi t}{i}-\frac{(\pi t)^3}{6i^3}+R_5\!\Bigl(\frac{\pi t}{i}\Bigr),
\qquad 
\Bigl|R_5\!\Bigl(\frac{\pi t}{i}\Bigr)\Bigr|\le \frac{(\pi |t|/i)^5}{120}.
\]
Factoring $\frac{\pi t}{i}$ from the denominator yields
\[
\sin\!\Bigl(\frac{\pi t}{i}\Bigr)
=\frac{\pi t}{i}\Bigl(1-\frac{\pi^2 t^2}{6i^2}+E_i(t)\Bigr),
\qquad 
|E_i(t)|\le \frac{(\pi |t|/i)^4}{120}\le \frac{\pi^4\varepsilon^4}{120}.
\]
Hence, for $|t|\le \varepsilon i$,
\[
\Bigl|\sin\!\Bigl(\frac{\pi t}{i}\Bigr)\Bigr|
\;\ge\;\frac{\pi |t|}{i}\Bigl(1-\frac{\pi^2\varepsilon^2}{6}-\frac{\pi^4\varepsilon^4}{120}\Bigr)
=\frac{\pi |t|}{i}\,c(\varepsilon).
\]
Division of the two Taylor series and inversion of the parenthesis by a Neumann series, justified by $c(\varepsilon)>0$ for $\varepsilon\in(0,\tfrac14]$, shows that the odd powers cancel and produces
\[
\frac{\sin(\pi x)}{\sin(\pi x/i)}
=(-1)^{m(i-1)}\,i\Bigl(1-\alpha_i(\pi t)^2\Bigr)+R_i^{(s)}(t),
\qquad \alpha_i=\frac{1}{6}\Bigl(1-\frac{1}{i^2}\Bigr),
\]
with the uniform remainder bound
\[
|R_i^{(s)}(t)|\le i\,C(\varepsilon)\,(\pi |t|)^4,
\qquad 
C(\varepsilon)=\frac{1}{120}\,c(\varepsilon)^{-2}.
\]
Squaring gives
\[
F(x,i)=\left(\frac{\sin(\pi x)}{\sin(\pi x/i)}\right)^{\!2}
= i^2\Bigl(1-2\alpha_i(\pi t)^2\Bigr)+R_i^{(F)}(t),
\]
with $|R_i^{(F)}(t)|\le i^2\,\widetilde C(\varepsilon)\,(\pi |t|)^4$ for some explicit $\widetilde C(\varepsilon)\asymp C(\varepsilon)$. 

For the resonance set, $|x/i-m|<\varepsilon$ is equivalent to $i\in I_m=(x/(m+\varepsilon),x/(m-\varepsilon))$, intersected with $(2,\infty)$ to respect the index range. Its length equals $|I_m|=\frac{2\varepsilon\,x}{m^2-\varepsilon^2}$. Therefore
\[
\sum_{m\ge \lceil \sqrt{x}\rceil}\! |I_m|
=\sum_{m\ge \lceil \sqrt{x}\rceil}\frac{2\varepsilon\,x}{m^2-\varepsilon^2}
\le \sum_{m\ge \lceil \sqrt{x}\rceil}\frac{2\varepsilon\,x}{(m-\varepsilon)^2}
\le \int_{\lceil \sqrt{x}\rceil-1}^{\infty}\frac{2\varepsilon\,x\,dt}{(t-\varepsilon)^2}
=\frac{2\varepsilon\,x}{\lceil \sqrt{x}\rceil-1-\varepsilon}.
\]
Since $(\lceil \sqrt{x}\rceil-1-\varepsilon)\ge \frac12\sqrt{x}$ for all $x\ge 1$ and $\varepsilon\le \tfrac14$, it follows that $\sum_{m\ge \lceil \sqrt{x}\rceil}|I_m|\le 4\varepsilon\sqrt{x}$.
\end{proof}
\end{lemma}

\begin{theorem}[RPF truncation with explicit remainder]\label{thm:RPF-trunc}
Let $i\ge 2$, $x\in\mathbb{R}$, $m:=\mathrm{round}(x/i)$, and $K\in\mathbb{N}_0$. Define the symmetric truncation
\[
S_{i,K}(x)\;:=\;\frac{1}{(x-i m)^{2}}+\sum_{r=1}^{K}\Biggl[\frac{1}{\bigl(x-i(m+r)\bigr)^{2}}+\frac{1}{\bigl(x-i(m-r)\bigr)^{2}}\Biggr].
\]
Then
\[
\left|\,F(x,i)-\frac{i^{2}}{\pi^{2}}\sin^{2}(\pi x)\,S_{i,K}(x)\,\right|
\ \le\ \frac{\sin^{2}(\pi x)}{\pi^{2}}\cdot \frac{2}{\,K+\tfrac12\,}.
\]
\end{theorem}

\begin{remark}[Integer arguments in the RPF truncation]
At integer $x=n$, the estimate is interpreted as $x\to n$. If $i\nmid n$, then $\sin(\pi n/i)\neq 0$ and no ambiguity occurs. If $i\mid n$, the removable $0/0$ singularity is resolved by the entire extension, yielding $F(n,i)=i^2$.
\end{remark}

\begin{proof}
By Proposition~\ref{prop:RPF}, $F(x,i)=(i^{2}/\pi^{2})\sin^{2}(\pi x)\sum_{k\in\mathbb{Z}}(x-ik)^{-2}$. For the chosen $m$, the distances satisfy $|x-i(m\pm r)|\ge (r-\tfrac12)i$ for $r\ge 1$. Hence the tail is bounded by
\[
\sum_{r>K}\frac{2}{\bigl(x-i(m\pm r)\bigr)^{2}}
\ \le\ \frac{2}{i^{2}}\sum_{r=K+1}^{\infty}\frac{1}{(r-\tfrac12)^{2}}
\ \le\ \frac{2}{i^{2}}\int_{K+\tfrac12}^{\infty}\frac{dr}{r^{2}}
\ =\ \frac{2}{i^{2}\,(K+\tfrac12)}.
\]
Multiplication by $(i^{2}/\pi^{2})\sin^{2}(\pi x)$ yields the claim.
\end{proof}

\begin{corollary}[Relative truncation]\label{cor:RPF-relative}
Let $x\notin\mathbb Z$, $i\ge 2$ and $K\in\mathbb N_0$, and set $m:=\mathrm{round}(x/i)$. Then
\[
\frac{\Bigl|\,F(x,i)-\frac{i^{2}}{\pi^{2}}\sin^{2}(\pi x)\,S_{i,K}(x)\,\Bigr|}{F(x,i)}
\ \le\ \frac{2}{K+\tfrac12}\cdot\frac{(x-i m)^2}{i^2}
\ \le\ \frac{1}{2(K+\tfrac12)}.
\]
The last inequality follows from $|x-i m|\le i/2$ for $m=\mathrm{round}(x/i)$. At integer arguments $x=n$ the inequality is to be interpreted in the removable-singularity sense, i.e. as the limit $x\to n$. In particular, for $K\ge 1$ the bound is at most $1/3$.
\end{corollary}

\noindent\textit{Algorithm for stable evaluation of $\mathcal{P}(x)$.} A reference implementation is provided in Appendix~\ref{app:code}.

\section{Smooth analogues of arithmetic functions: setup}

\subsection{Motivation}
The prime-indicating property of $\mathcal{P}(x)$ comes with limited smoothness.
The discontinuities in its second derivative are a direct result of the sharp cutoff in the summation limit $\lceil\sqrt{x}\rceil$.
To obtain globally smooth ($C^\infty$) functions, the sharp cutoff is replaced by a smooth transition, yielding a convergent infinite series.
The representation in \eqref{eq:Px} extends to weighted Fejér sums that reproduce divisor sums at integer arguments.
The value of $\mathcal{P}(n)$ for a composite integer $n$ is $\frac{1}{n}\sum d^2$, which is not immediately intuitive.
This motivates a generalization to create smooth functions that correspond to more classical arithmetic functions, such as the divisor-counting function $\tau(n)$ (the number of positive divisors of $n$) or the sum-of-divisors function $\sigma(n)$ (the sum of positive divisors of $n$).
Definitions of $\mathcal{S}_W$ and $\mathcal{P}_W$ are given below.

\subsection{Structure of the construction}
The weighted core sum is defined as
\[ \mathcal{S}_W(x; \kappa) = \sum_{i=2}^{\infty} \phi_{\kappa}\!\left(\frac{i}{x+1}\right) \, W(i) \, F(x,i), \]
where $F(x,i)$ is the Fejér term from \eqref{eq:Fxi}, $W(i)$ is a weighting function, and $\phi_{\kappa}$ is a smooth cutoff. The quantity
\[ \mathcal{P}_W(x; \kappa) := \mathcal{S}_W(x; \kappa) - N_W(x) \]
uses $N_W$ as the normalization that enforces the desired vanishing at primes on integer input. In particular, for $W(i)=i^{-2}$ one sets $N_W(x)\equiv 1$, which yields the $\tau$-analogue $\mathcal P_\tau=\mathcal S_{i^{-2}}-1$, and for $W(i)=i^{-1}$ one sets $N_W(x)=x$, which yields the $\sigma$-analogue $\mathcal P_\sigma=\mathcal S_{i^{-1}}-x$. For integers $n$, these choices reproduce the targets $\tau(n)-2$ and $\sigma(n)-n-1$ in the steep-cutoff limit $\kappa\to\infty$.

\subsection{Normalization}
The factor $1/x$ present in the definition of $\mathcal{P}(x)$ is omitted in the smooth setting. This omission facilitates the construction of analogues for classical, unnormalized arithmetic functions. At integer arguments $n$, the identity $F(n,d)=d^2$ turns the core sum $\sum_d W(d)F(n,d)$ into a divisor sum whose contributions are directly controlled by the weighting function $W$.

For instance, $W(d)=d^{-2}$ yields a unit contribution from each divisor (targeting $\tau(n)$), while $W(d)=d^{-1}$ yields a contribution of $d$ from each divisor (targeting $\sigma(n)$). The resulting prime indicator $\mathcal{P}_\sigma(x) = \mathcal{S}_\sigma(x) - x$ uses an additive comparison rather than a multiplicative one (e.g., $\mathcal{S}_\sigma(x)/x - 1$). This choice directly models the vanishing criterion for primes at the integer level, $\sigma(n)-n-1=0$, without introducing a denominator that might alter the analytic structure away from the integers.

\begin{remark}[Integer evaluation and connection to divisor sums]\label{rem:integer-eval-smooth}
At integer arguments $n$, the filter property $F(n,i)=\mathbf{1}_{i \mid n} \cdot i^2$ simplifies the core sum $\mathcal{S}_W(n;\kappa)$ to
\[
\mathcal{S}_W(n;\kappa) = \sum_{\substack{d \mid n \\ d \ge 2}} \phi_{\kappa}\left(\frac{d}{n+1}\right) W(d) d^2.
\]
The choice of weights is thus designed to reproduce classical arithmetic functions in the steep-cutoff limit $\kappa\to\infty$, where $\phi_\kappa(d/(n+1))\to 1$ for any divisor $d\le n$:
\begin{itemize}
    \item With $W(i)=i^{-2}$ (for $\mathcal{P}_\tau$), the limit of the sum becomes $\sum_{d|n, d\ge2} 1 = \tau(n)-1$.
    \item With $W(i)=i^{-1}$ (for $\mathcal{P}_\sigma$), the limit of the sum becomes $\sum_{d|n, d\ge2} d = \sigma(n)-1$.
\end{itemize}
For any finite $\kappa>0$ and any odd prime $p$, the only contributing divisor is $d=p$. The core sums evaluate to $\phi_{\kappa}(p/(p+1))$ for the $\tau$-case and $p\,\phi_{\kappa}(p/(p+1))$ for the $\sigma$-case. Since $\phi_\kappa(u)<1$ for $u<1$, the resulting indicator values are strictly negative. Exact vanishing at primes is thus a property of the sharp-cutoff limit, not of the smooth functions for finite $\kappa$.
\end{remark}

\subsection{A Smooth Cutoff Function}
A cutoff function is required that is essentially $1$ for $u \le 1$ and decays rapidly for $u > 1$.
The choice $u=i/(x+1)$ includes indices with $i\lesssim x$ and, together with $\phi_\kappa$, yields absolute convergence of the infinite series.
A function based on the hyperbolic tangent has these properties and is of class $C^\infty$.
\begin{definition}[Smooth Transition Cutoff Function]
Let $\kappa>0$ be a steepness parameter. The modified cutoff function $\phi_{\kappa}\colon\R\to(0,1)$ is defined as
\begin{equation}\label{eq:phimod}
\phi_{\kappa}(u) = \frac{1 - \tanh\big(\kappa (u - 1)\big)}{2}.
\end{equation}
\end{definition}

\begin{lemma}[Cutoff limit and value at the threshold]
For every $\kappa>0$, $\phi_\kappa(1)=\tfrac12$. For each fixed $\delta\in(0,1)$, the convergence $\phi_\kappa(u)\uparrow 1$ as $\kappa\to\infty$ is uniform on $(-\infty,1-\delta]$, and $\phi_\kappa(u)\downarrow 0$ as $\kappa\to\infty$ is uniform on $[1+\delta,\infty)$. Pointwise, $\phi_\kappa(u)\to \mathbf 1_{\{u<1\}}$ for all $u\neq 1$.
\end{lemma}

There is a smooth transition from $1$ to $0$ centered at $u=1$.
The logistic form $\phi_\kappa(u)=\bigl(1+e^{2\kappa(u-1)}\bigr)^{-1}$ implies the monotone limit
$\phi_\kappa(u)\uparrow 1$ as $\kappa\to\infty$ for every $u<1$, with the quantitative bound
\[
1-\phi_\kappa(u)=\frac{e^{2\kappa(u-1)}}{1+e^{2\kappa(u-1)}}\le e^{2\kappa(u-1)}.
\]
The cutoff is applied to $u=i/(x+1)$. For an integer $x=n$ and any divisor $d\mid n$, one has
$u=d/(n+1)\le n/(n+1)=1-\frac{1}{n+1}$, hence
\[
1-\phi_\kappa\!\Bigl(\frac{d}{n+1}\Bigr)\ \le\ e^{-\,\frac{2\kappa}{\,n+1\,}},
\]
uniformly over all divisors $d$ (including $d=n$). Using $x+1$ instead of $x$ keeps $u$ uniformly away from the threshold at $1$ for $x>0$ and avoids degenerate cases.

\begin{remark}[On the choice of cutoff]\label{rem:cutoff-choice}
The dependence $\phi_\kappa(i/(x+1))$ cannot be replaced by a purely $x$-dependent weight without losing absolute convergence and the correct integer limits. In particular, using a factor $\psi_\kappa(x)$ independent of $i$ would turn the $\tau$-series into a sum with generic term $F(x,i)/i^2\le 1$, which is not summable, and the $\sigma$-series into a sum with terms $\asymp i$, which diverges even faster. The variable $u=i/(x+1)$ is essential to both analytic and arithmetic control.
\end{remark}

\begin{figure}[!htbp]
\centering
\includegraphics[width=0.8\textwidth]{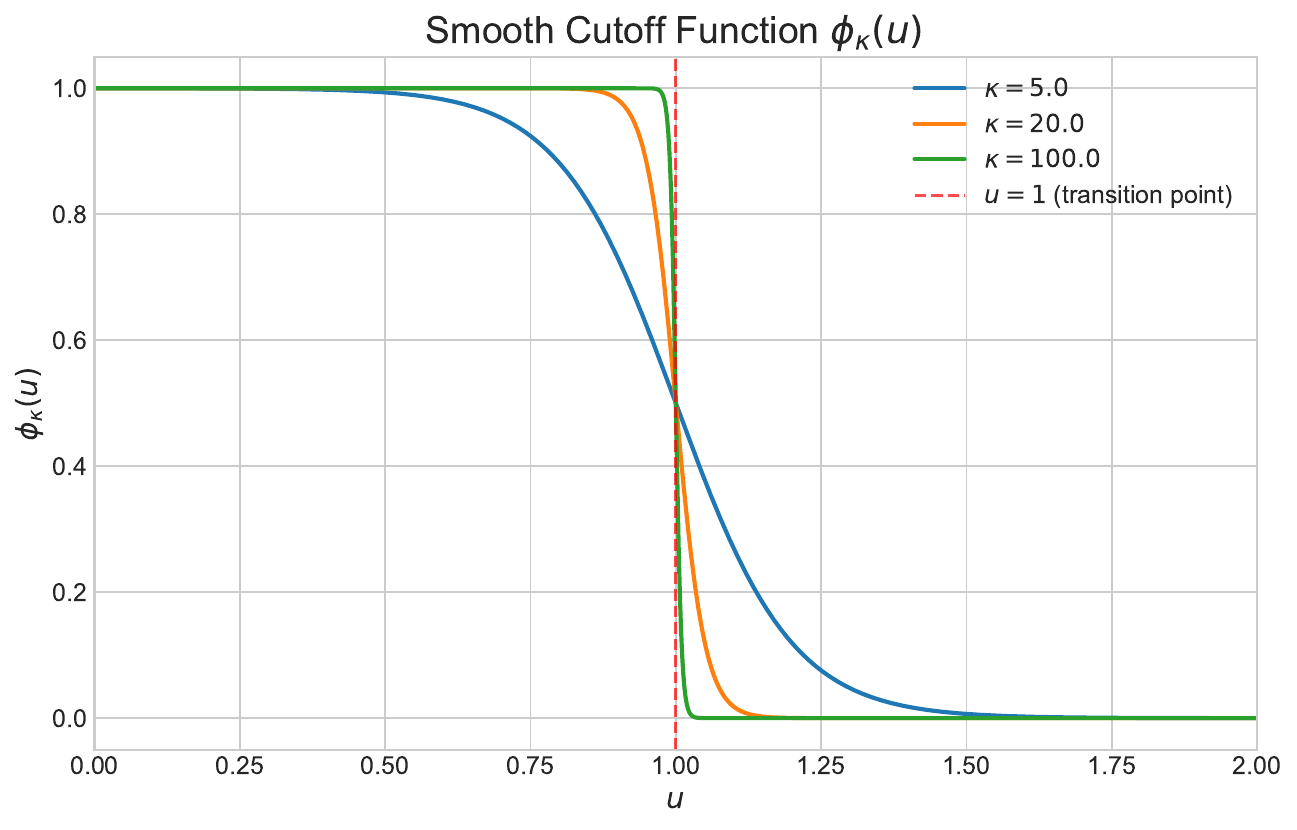}
\caption{The smooth cutoff function $\phi_{\kappa}(u)$ for different values of the steepness parameter $\kappa$.
A smooth transition from 1 to 0 occurs around $u=1$.
For larger $\kappa$, the transition becomes steeper while $C^\infty$-smoothness is preserved.}
\label{fig:cutoff}
\end{figure}

\FloatBarrier

\subsection{Numerical complexity and truncation with a smooth cutoff}

\begin{theorem}[Uniform and pointwise truncation bounds with explicit constants]\label{thm:truncation}
Let $K=[a,b]\subset(0,\infty)$ be compact and $\kappa>0$. Set $c_K:=\frac{2\kappa}{\,b+1\,}$ and $r:=e^{-c_K}\in(0,1)$. Then for every integer $M\in\mathbb{N}$ with $M\ge b+1$,
\begin{align*}
\sup_{x\in K}\ \sum_{i>M}\phi_\kappa\!\left(\frac{i}{x+1}\right)\frac{|F(x,i)|}{i^2}
&\ \le\ \frac{r^{\,M+1}}{1-r},\\[4pt]
\sup_{x\in K}\ \sum_{i>M}\phi_\kappa\!\left(\frac{i}{x+1}\right)\frac{|F(x,i)|}{i}
&\ \le\ r^{\,M+1}\!\left(\frac{M+1}{1-r}+\frac{r}{(1-r)^2}\right)
\ \le\ \frac{(M+2)\,r^{\,M+1}}{(1-r)^2}.
\end{align*}
For any $\varepsilon\in(0,1)$ a sufficient uniform choice for the $\tau$-case is
\[
M_\tau^{\mathrm{unif}}=\left\lceil -\,1+\frac{1}{c_K}\log\!\Bigl(\frac{1}{(1-e^{-c_K})\,\varepsilon}\Bigr)\right\rceil.
\]
For the $\sigma$-case one may select the minimal integer $M\ge b+1$ such that
\[
\frac{(M+2)\,e^{-c_K(M+1)}}{(1-e^{-c_K})^2}\ \le\ \varepsilon,
\]
which exists because the left-hand side tends to $0$ as $M\to\infty$ and is eventually strictly decreasing in $M$. For fixed $x>0$, the pointwise analogues follow upon replacing $c_K$ by $c_x:=2\kappa/(x+1)$ and choosing any integer $M\ge \lfloor x\rfloor+1$ so that $i/(x+1)\ge 1$ holds for all $i>M$.
\begin{proof}
If $M\ge b+1$ then $i>M$ implies $u=i/(x+1)\ge 1$ for all $x\in[a,b]$. Using the exact logistic form,
\[
\phi_\kappa(u)=\frac{1}{1+e^{2\kappa(u-1)}}\ \le\ e^{-2\kappa(u-1)}\ \le\ e^{-c_K i}=r^{\,i}.
\]
With $|F(x,i)|\le i^2$ and $\sum_{i>M}r^i=r^{M+1}/(1-r)$ the first bound follows. For the second bound, $\sum_{i>M} i\,r^i=r^{M+1}\bigl(\frac{M+1}{1-r}+\frac{r}{(1-r)^2}\bigr)\le \frac{(M+2)\,r^{M+1}}{(1-r)^2}$. The remaining statements are immediate.
\end{proof}
\end{theorem}

\begin{remark}[On the tail constants]
The bounds in Theorem~\ref{thm:truncation} already use the exact logistic inequality $\phi_\kappa(u)\le e^{-2\kappa(u-1)}$ for $u\ge 1$. No additional prefactors are required, and any re-centering of the tail only improves the global constant without changing the geometric rate.
\end{remark}

Consequently, the cost per evaluation at a given $x$ is $O\!\bigl(x+\tfrac{x}{\kappa}\log(1/\varepsilon)\bigr)$ for $\mathcal P_\tau$ and $O\!\bigl(x+\tfrac{x}{\kappa}[\log(1/\varepsilon)+\log(1+x)]\bigr)$ for $\mathcal P_\sigma$. Near resonances $x/i\approx\mathbb{Z}$, numerical stability is improved by switching to the cosine-polynomial representation of $F(x,i)$, as in the baseline pseudocode.

\section{Application 1: A Smooth Divisor-Counting Function (\texorpdfstring{$\mathcal{P}_{\tau}$}{P_tau})}

\subsection{Motivation and Construction}
A smooth analogue of the divisor-counting function $\tau(n)$ is constructed. For an integer $n$, the Fejér term $F(n,d)$ evaluates to $d^2$ if $d$ divides $n$. A unit contribution per divisor is obtained by scaling with $d^{-2}$, which leads to the weighting function $W(i)=1/i^2$. The resulting sum over $i \ge 2$ approximates $\tau(n)-1$ (the number of divisors excluding $1$). A subtraction by $1$ is then applied so that the expression vanishes at primes (where $\tau(n)-1=1$).

\begin{definition}[Smooth Divisor-Counting Function]
For $x>0$ and a steepness parameter $\kappa>0$, define the smooth divisor-counting analogue
\begin{equation}\label{eq:Ptau}
\mathcal{P}_{\tau}(x; \kappa) = \left( \sum_{i=2}^{\infty} \phi_{\kappa}\!\left(\frac{i}{x+1}\right) \cdot \frac{F(x,i)}{i^2} \right) - 1,
\end{equation}
where $\phi_{\kappa}$ is given by \eqref{eq:phimod} and $F(x,i)$ by \eqref{eq:Fxi}.
\end{definition}

\begin{lemma}[Explicit cutoff-derivative bounds]\label{lem:cutoff-derivatives}
Let $\kappa>0$ and $K=[a,b]\subset(0,\infty)$ be compact. For every integer $r\ge 0$ and all $i\ge 2$,
\[
\sup_{x\in K}\,\Bigl|\partial_x^r\,\phi_{\kappa}\!\Bigl(\frac{i}{x+1}\Bigr)\Bigr|
\ \le\ \frac{(2\kappa)^r\,r!\,B_r}{(a+1)^{r+1}}\; i^{\,r}\, \exp\!\Bigl(-\frac{2\kappa}{\,b+1\,}\, i\Bigr),
\]
where $B_r>0$ depends only on $r$. Equivalently, one may set
\[
C_{r,K,\kappa}:=\frac{(2\kappa)^r\,r!\,B_r}{(a+1)^{r+1}},\qquad c_{K,\kappa}:=\frac{2\kappa}{\,b+1\,},
\]
so that $\sup_{x\in K}|\partial_x^r\,\phi_{\kappa}(i/(x+1))|\le C_{r,K,\kappa}\, i^{r}\, e^{-c_{K,\kappa}\, i}$.
\begin{proof}
Write $\phi_\kappa(u)=(1+e^{2\kappa(u-1)})^{-1}$ and $u(x)=i/(x+1)$. For $x\in K$ one has $u(x)\ge i/(b+1)$, hence $\phi_\kappa(u(x))\le e^{-2\kappa(u(x)-1)}\le \exp(-(2\kappa/(b+1))\,i)$. Set $y(x)=\kappa(u(x)-1)$ and use Faà di Bruno’s formula for derivatives of $\tanh$ with $\phi_\kappa=(1-\tanh y)/2$. Since $u^{(j)}(x)=(-1)^j j!\,i\,(x+1)^{-(j+1)}$, it follows that $\sup_{x\in K}|u^{(j)}(x)|\le j!\,i\,(a+1)^{-(j+1)}$. Moreover, $\mathrm{sech}^2 y(x)\le 4e^{-2|y(x)|}\le 4\,e^{-2\kappa(u(x)-1)}$. Collecting factors gives
\[
\sup_{x\in K}\bigl|\partial_x^r\,\tanh(y(x))\bigr|\ \le\ \frac{(2\kappa)^r\,r!\,B_r}{(a+1)^{r+1}}\,i^{\,r}\,\exp\!\Bigl(-\frac{2\kappa}{\,b+1\,}\,i\Bigr),
\]
and the same upper bound applies to $\partial_x^r\phi_\kappa(u(x))$ up to an absolute factor absorbed into $B_r$. This proves the stated bound.
\end{proof}
\end{lemma}

\subsection{Properties of $\mathcal{P}_{\tau}$}

\begin{proposition}[Fejér derivative bounds]\label{prop:fejer-bounds}
For $i\ge2$ and $x\in\R$, let
\[
F(x,i)= i + 2\sum_{k=1}^{i-1}(i-k)\cos\!\left(\frac{2\pi k x}{i}\right).
\]
Then for every integer $l\ge0$ there exist constants $A_l,B_l>0$, depending only on $l$, such that
\[
\sup_{x\in\R}\bigl| \partial_x^l F(x,i) \bigr|
\ \le\ \frac{2}{(l+1)(l+2)}(2\pi)^l\,i^{2}\;+\;B_l(2\pi)^l\,i,
\qquad \forall i\ge2.
\]
In particular, for $l=0$ one has $|F(x,i)|\le i^2$, and
\[
\sup_{x\in\R}\bigl| \partial_x^l (F(x,i)/i^2) \bigr| \ \le\ \frac{2}{(l+1)(l+2)}(2\pi)^l\;+\;B_l(2\pi)^l\,i^{-1},
\]
\[
\sup_{x\in\R}\bigl| \partial_x^l (F(x,i)/i) \bigr| \ \le\ \frac{2}{(l+1)(l+2)}(2\pi)^l\,i\;+\;B_l(2\pi)^l.
\]
\begin{proof}
Differentiation yields
\[
\partial_x^l F(x,i)=2\sum_{k=1}^{i-1}(i-k)\left(\frac{2\pi k}{i}\right)^{\!l}
\!\times\!
\begin{cases}
\cos(2\pi k x/i), & l\ \text{even},\\
-\sin(2\pi k x/i), & l\ \text{odd}.
\end{cases}
\]
Since $|\cos(\cdot)|\le 1$ and $|\sin(\cdot)|\le 1$,
\[
\bigl|\partial_x^l F(x,i)\bigr|
\le 2\left(\frac{2\pi}{i}\right)^{\!l}\sum_{k=1}^{i-1}(i-k)k^l.
\]
Faulhaber's formula gives $\sum_{k=1}^{i-1}(i-k)k^l=\frac{i^{l+2}}{(l+1)(l+2)}+R_l(i)$ with $|R_l(i)|\le A_l\,i^{l+1}$. Consequently,
\[
\sup_{x\in\R}\bigl|\partial_x^l F(x,i)\bigr|
\le \frac{2}{(l+1)(l+2)}(2\pi)^l\,i^{2}+B_l(2\pi)^l\,i,
\]
where one may take $B_l\ge A_l$ independent of $i$. The stated corollaries follow by dividing by $i^2$ or $i$.
\end{proof}
\end{proposition}

\begin{remark}[Sharper leading constant]
Faulhaber’s formula implies, for all fixed $l\ge 0$ and integers $i\ge 2$,
\[
\sum_{k=1}^{i-1}(i-k)k^l=\frac{i^{l+2}}{(l+1)(l+2)}+R_{l}(i),\qquad |R_{l}(i)|\le A_l\,i^{l+1},
\]
for some $A_l>0$. Consequently,
\[
\sup_{x\in\mathbb{R}}|\partial_x^l F(x,i)|
\le \frac{2}{(l+1)(l+2)}(2\pi)^l\, i^{2}\;+\;B_l\,(2\pi)^l\, i.
\]
Thus the leading constant equals $2/[(l+1)(l+2)]$ with a linear remainder term.
\end{remark}

\begin{remark}[Uniform $M$-test blueprint]\label{rem:mtest-blueprint}
On the real axis, $|F(x,i)|\le i^2$ and $|\partial_x^l F(x,i)|\le C_l(2\pi)^l i^2$. For the cutoff (here and below $r$ denotes the order of differentiation),
\[
\sup_{x\in K}\bigl|\partial_x^r \phi_\kappa(i/(x{+}1))\bigr|\ \ll_{r,K,\kappa}\ i^{\,r}\,e^{-c i}
\]
holds on every compact $K\subset(0,\infty)$ for some $c=c_{K,\kappa}>0$. By Leibniz' rule,
\[
\sup_{x\in K}\bigl|\partial_x^r\bigl(\phi_\kappa(i/(x{+}1))\,F(x,i)/i^\beta\bigr)\bigr|\ \ll_{r,K,\kappa}\ i^{\,r+\alpha}\,e^{-c i},
\]
with $\beta\in\{1,2\}$ and $\alpha\in\{0,1\}$, hence the series of derivatives is summable in $i$. The Weierstrass $M$-test yields local uniform convergence of all derivative series and thus $C^\infty$-regularity.
\end{remark}

\begin{lemma}[Local uniform convergence implies smoothness]\label{lem:local_uniform_smooth}
Let $(h_i)_{i\ge2}$ be a sequence of $C^\infty$ functions on $(0,\infty)$. If for every compact $K\subset(0,\infty)$ and every $j\ge0$ the series $\sum_{i\ge2}\partial_x^j h_i$ converges uniformly on $K$, then $H=\sum_{i\ge2} h_i$ defines a $C^\infty$ function on $(0,\infty)$ and derivatives may be taken termwise on compacts.
\end{lemma}

\begin{proposition}[Properties of $\mathcal{P}_{\tau}$]
The function $\mathcal{P}_{\tau}(x;\kappa)$ is of class $C^\infty$ on $(0,\infty)$. For any integer $n\ge2$,
\[
\lim_{\kappa\to\infty} \mathcal{P}_{\tau}(n; \kappa) = \tau(n)-2,
\]
which vanishes if and only if $n$ is prime.
\end{proposition}
\begin{proof}
To establish $C^\infty$-regularity on $(0,\infty)$, it suffices to show that for any compact set $K\subset(0,\infty)$ and any integer $j\ge0$, the series of $j$-th derivatives converges uniformly on $K$. Let $h_i(x) = \phi_{\kappa}(i/(x+1)) F(x,i)/i^2$. By the Leibniz rule,
\[
\partial_x^j h_i(x) = \sum_{q=0}^j \binom{j}{q} \left( \partial_x^q \phi_{\kappa}\!\left(\frac{i}{x+1}\right) \right) \left( \partial_x^{j-q} \frac{F(x,i)}{i^2} \right).
\]
From Lemma~\ref{lem:cutoff-derivatives}, $|\partial_x^q \phi_{\kappa}(i/(x+1))| \le C_{q,K} \, i^q \, e^{-c_K i}$ for constants $C_{q,K}, c_K > 0$. From Proposition~\ref{prop:fejer-bounds}, the derivative $|\partial_x^{j-q} (F(x,i)/i^2)|$ is bounded by a constant, say $D_{j-q}$. Combining these gives a uniform majorant on $K$:
\[
\sup_{x\in K} |\partial_x^j h_i(x)| \le \sum_{q=0}^j \binom{j}{q} (C_{q,K} \, i^q \, e^{-c_K i}) D_{j-q} \le M_{j,K} \, i^{j} \, e^{-c_K i}.
\]
This majorant is summable over $i$ for any fixed $j$. By the Weierstrass M-test, each derivative series converges uniformly on $K$, which implies that $\mathcal{P}_{\tau}$ is a $C^\infty$ function.

For the integer limit, if $n\in\mathbb Z$, Proposition~\ref{prop:fejer-properties} gives $F(n,i)/i^2 = \mathbf{1}_{i\mid n}$. The series becomes a finite sum:
\[
\mathcal{P}_{\tau}(n;\kappa)=\sum_{\substack{d\mid n\\ d\ge2}} \phi_{\kappa}\!\Bigl(\frac{d}{n+1}\Bigr)-1.
\]
As $\kappa\to\infty$, $\phi_{\kappa}(u)\to\mathbf 1_{u<1}$. Since $d/(n+1)<1$ for any divisor $d$ of $n$, the limit is $\#\{d\mid n: d\ge2\}-1 = \tau(n)-2$.
\end{proof}

\begin{definition}[Zero terminology]\label{def:zero-terminology}
For a prime $p$, a \emph{companion zero at $p$} is any real zero $x\neq p$ of a smooth prime indicator with $|x-p|<1$. A family of zeros $x_p^{(\pm)}(\kappa)$ parameterized by the steepness $\kappa$ is called a \emph{prime-zero approximant} if $x_p^{(\pm)}(\kappa)\to p$ as $\kappa\to\infty$.
In the present setting, $x=p$ is not a zero of the smooth indicators $\mathcal P_\tau(\cdot;\kappa)$ or $\mathcal P_\sigma(\cdot;\kappa)$ for any finite $\kappa$; the integer zero at $p$ belongs only to the sharp-cutoff indicator $\mathcal P$.
Conjecture~\ref{conj:companion-zeros} concerns two approximants flanking $p$ for $\mathcal P_\tau$, whereas Conjecture~\ref{conj:Psigma-companions} states an asymmetric pair for $\mathcal P_\sigma$ (left approximant; right companion at a $\kappa$-independent distance). The terminology is used consistently below without further repetition.
\end{definition}

\begin{conjecture}[Companion zeros near odd primes]\label{conj:companion-zeros}
For each odd prime $p\ge 3$, it is conjectured that there exist $\kappa_0(p)\ge 1$ and $\varepsilon_p>0$ such that, for all $\kappa\ge \kappa_0(p)$, the function $\mathcal P_\tau(\cdot;\kappa)$ admits two simple real zeros, one in $(p-\varepsilon_p,p)$ and one in $(p,p+\varepsilon_p)$. Numerical computations suggest the distances to $p$ are $O\bigl(\mathrm e^{-\kappa/(p+1)}\bigr)$ as $\kappa\to\infty$.
\end{conjecture}

\begin{figure}[!htbp]
\centering
\includegraphics[width=0.9\textwidth]{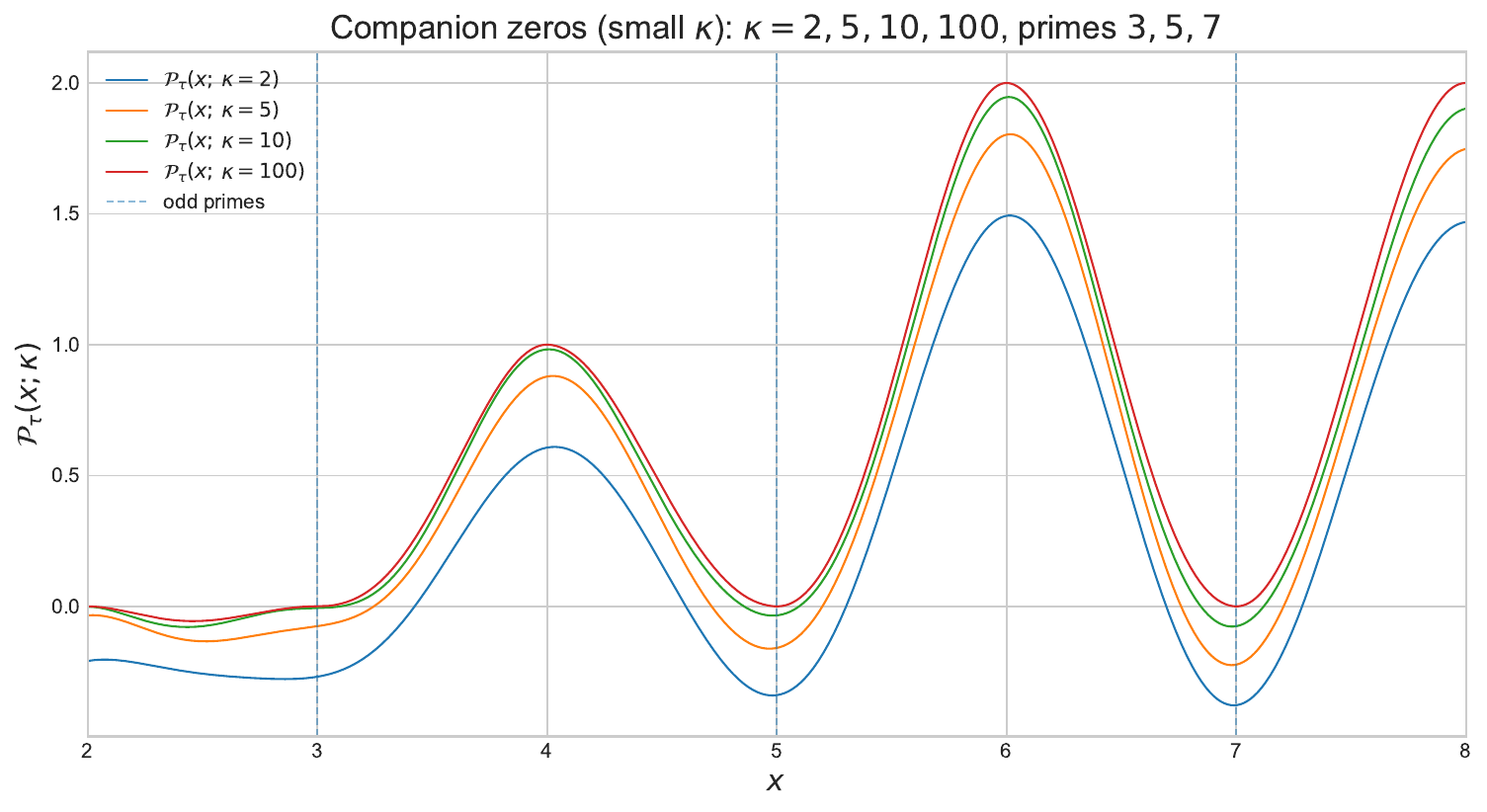}
\caption{Companion zeros of the smooth divisor–counting analogue $\mathcal{P}_{\tau}(x;\kappa)$
on $[2,8]$ for $\kappa\in\{2,5,10,100\}$. Vertical dashed lines mark the odd primes $3,5,7$.
Consistent with Conjecture~\ref{conj:companion-zeros}, for larger $\kappa$ a pair of real zeros
appears on either side of each odd prime and moves toward $p$; the value at $x=p$
is strictly negative for finite $\kappa$ and decays like $\mathrm e^{-2\kappa/(p+1)}$.}
\label{fig:PtauCompanions}
\end{figure}

As a numerical illustration of Conjecture~\ref{conj:companion-zeros}, 
Figure~\ref{fig:PtauCompanions} plots $\mathcal P_\tau(x;\kappa)$ on $[2,8]$ for 
$\kappa\in\{2,5,10,100\}$; the pair of real zeros flanking each odd prime emerges and tightens 
with~$\kappa$, while the central value at $x=p$ remains negative for finite~$\kappa$.


\begin{lemma}[Local quadratic lower bound for $\mathcal P_\tau$ at an odd prime]\label{lem:Ptau-local-lb}
Let $p\ge 3$ be prime. There exist $\kappa_0(p)\ge 1$, $\varepsilon_p>0$, and $c_p>0$ such that, for all $\kappa\ge \kappa_0(p)$ and $|x-p|\le \varepsilon_p$,
\[
\mathcal P_\tau(x;\kappa)\;=\;-\delta_p^{(\tau)}(\kappa)\;+\;B_p(\kappa)\,(x-p)^2\;+\;O\bigl((x-p)^3\bigr),
\qquad \delta_p^{(\tau)}(\kappa):=1-\phi_\kappa\!\Bigl(\tfrac{p}{p+1}\Bigr)\asymp \mathrm e^{-2\kappa/(p+1)}.
\]
As $\kappa\to\infty$,
\[
B_p(\kappa)\ \longrightarrow\ B_p^{(\infty)}
=\sum_{\substack{2\le i\le p\\ i\neq p}}\frac{\pi^2}{i^2\sin^2(\pi p/i)}
+\frac{1}{2}\,\frac{\pi^2}{(p+1)^2\sin^2\!\bigl(\pi/(p+1)\bigr)}
-\frac{\pi^2}{3}\Bigl(1-\frac{1}{p^2}\Bigr),
\]
and $B_p^{(\infty)}>0$. In particular,
\[
B_p^{(\infty)} \;\ge\; \frac{\pi^2}{4} \;+\; \sum_{i=3}^{p-1}\frac{\pi^2}{i^2} \;+\; \frac12 \;-\; \frac{\pi^2}{3}\Bigl(1-\frac{1}{p^2}\Bigr),
\]
and for $p=3$ this bound reduces to $ \frac{7\pi^2}{432}>0$.
\end{lemma}

\begin{proof}
Set $t=x-p$. Since each $F(\cdot,i)$ and $\phi_\kappa$ is $C^\infty$, a second-order Taylor expansion at $x=p$ is legitimate and uniform on $2\le i\le p+1$.

\emph{(i) $2\le i\le p-1$.} Because $\sin(\pi p/i)\neq 0$,
\[
F(x,i)=\frac{\sin^2(\pi x)}{\sin^2(\pi x/i)}=\frac{\pi^2}{\sin^2(\pi p/i)}\,t^2+O(t^3),
\]
and therefore
\[
\phi_\kappa\!\Bigl(\frac{i}{x+1}\Bigr)\frac{F(x,i)}{i^2}
=\Bigl(\phi_\kappa\!\Bigl(\tfrac{i}{p+1}\Bigr)+O(t)\Bigr)\Bigl(\frac{\pi^2}{i^2\sin^2(\pi p/i)}\,t^2+O(t^3)\Bigr).
\]
Since $\phi_\kappa(\tfrac{i}{p+1})\to 1$ as $\kappa\to\infty$, the quadratic coefficient tends to $\pi^2/(i^2\sin^2(\pi p/i))$.

\emph{(ii) $i=p$.} With $x=p+t$,
\[
\frac{\sin(\pi x)}{\sin(\pi x/p)}=p\Bigl(1-\alpha_p(\pi t)^2\Bigr)+O(t^4),\qquad
\alpha_p=\frac{1}{6}\Bigl(1-\frac{1}{p^2}\Bigr),
\]
hence
\[
\frac{F(x,p)}{p^2}=1-2\alpha_p(\pi t)^2+O(t^4),
\]
and the quadratic coefficient equals $-\,\dfrac{\pi^2}{3}\Bigl(1-\dfrac{1}{p^2}\Bigr)$.

\emph{(iii) $i=p+1$.} Here $u(x)=\tfrac{p+1}{x+1}=1-\tfrac{t}{p+1}+O(t^2)$ and
\[
\phi_\kappa\!\Bigl(\frac{p+1}{x+1}\Bigr)=\frac12+O(t) \quad (\text{for fixed }\kappa).
\]
Since
\[
\frac{F(x,p+1)}{(p+1)^2}=\frac{\pi^2}{(p+1)^2\sin^2(\pi/(p+1))}\,t^2+O(t^3),
\]
the product has quadratic coefficient \emph{exactly}
\[
\frac{1}{2}\cdot\frac{\pi^2}{(p+1)^2\sin^2(\pi/(p+1))},
\]
because the $O(t)$ part of $\phi_\kappa$ multiplies a function without linear term.

Summing (i)–(iii) and subtracting the constant $1$ in the definition of $\mathcal P_\tau$ yields $B_p^{(\infty)}$ as stated. The lower bound follows from $\sin^2\theta\le 1$ and $\sin(\pi p/2)=\pm 1$ for odd $p$. For $p=3$,
\(
\frac{\pi^2}{4}+\frac{\pi^2}{16}-\frac{8\pi^2}{27}=\frac{7\pi^2}{432}>0.
\)
\end{proof}

\begin{corollary}[Quantitative flank zeros near $p$]
For $p\ge3$ and $\kappa\ge \kappa_0(p)$, there exist two simple real zeros of $\mathcal P_\tau(\cdot;\kappa)$ in a neighborhood of $p$, with distances
\[
|x-p|\;=\;\sqrt{\frac{\delta_p^{(\tau)}(\kappa)}{B_p(\kappa)}}\;+\;O\!\left(\frac{\delta_p^{(\tau)}(\kappa)}{B_p(\kappa)}\right)
\;\asymp\; \mathrm e^{-\kappa/(p+1)}.
\]
\end{corollary}

\begin{figure}[!htbp]
\centering
\includegraphics[width=0.9\textwidth]{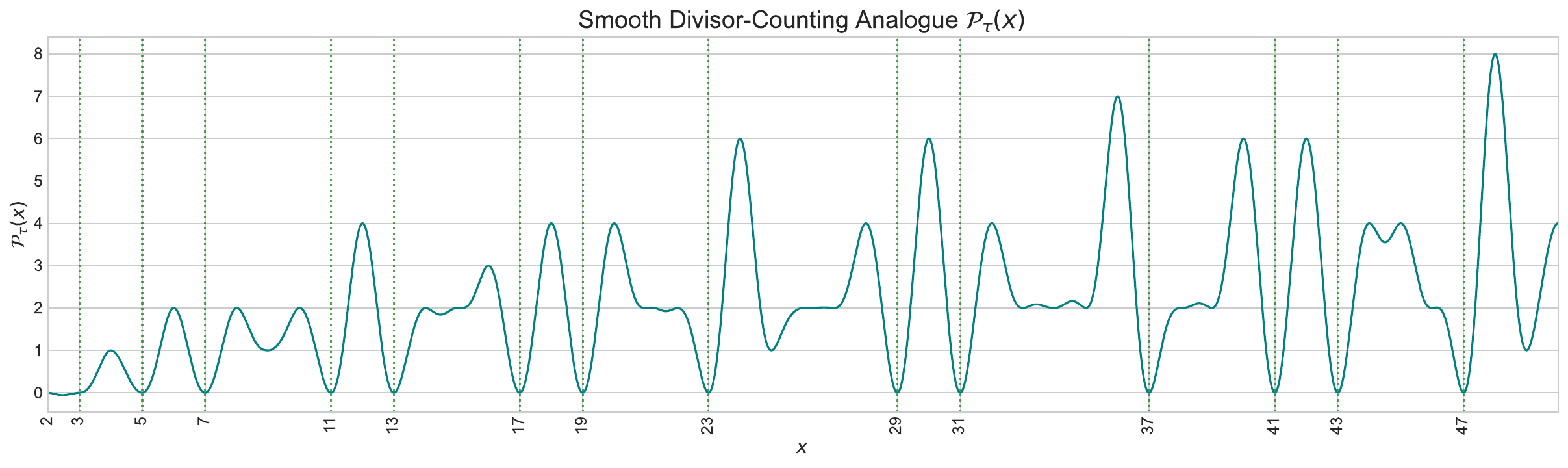}
\caption{Numerical profile of $\mathcal{P}_{\tau}(x; \kappa)$ (default $\kappa=1000$ unless stated otherwise).
The function is $C^\infty$ and approaches $0$ at all integer primes as $\kappa\to\infty$. The local shape near odd primes is consistent with Conjecture~\ref{conj:companion-zeros}.}
\label{fig:P_tau}
\end{figure}

\FloatBarrier

\section{Application 2: A Smooth Divisor-Sum Function (\texorpdfstring{$\mathcal{P}_{\sigma}$}{P_sigma})}

\subsection{Motivation and Construction}
A smooth analogue of the sum-of-divisors function, $\sigma(n)$, can be constructed.
To achieve this, it is required that each divisor $d$ contributes its own value $d$ to the sum.
Since $F(n,d)=d^2$, the natural choice for the weighting function is $W(i)=1/i$.
\begin{definition}[Smooth Divisor-Sum Function]
For $x>0$ and a steepness parameter $\kappa>0$, define
\begin{equation}\label{eq:Psigma}
\mathcal{P}_{\sigma}(x; \kappa) = \left( \sum_{i=2}^{\infty} \phi_{\kappa}\!\left(\frac{i}{x+1}\right) \cdot \frac{F(x,i)}{i} \right) - x.
\end{equation}
\end{definition}
For integer arguments $n$, the summation part evaluates to the sum of all divisors except 1, which is $\sigma(n)-1$. Subtraction of $n$ yields an expression that vanishes precisely at primes, for which $\sigma(n)-1=n$.
\subsection{Properties of $\mathcal{P}_{\sigma}$}

\begin{proposition}[Properties of $\mathcal{P}_{\sigma}$]
The function $\mathcal{P}_{\sigma}(x;\kappa)$ is $C^\infty$ on $(0,\infty)$. For any integer $n\ge2$,
\[
\lim_{\kappa\to\infty} \mathcal{P}_{\sigma}(n; \kappa) = \sigma(n) - n - 1,
\]
and the right-hand side vanishes if and only if $n$ is prime.
\end{proposition}

\begin{proof}
The proof of $C^\infty$-regularity on $(0,\infty)$ follows the same argument as for $\mathcal{P}_{\tau}$. Let $h_i(x) = \phi_{\kappa}(i/(x+1)) F(x,i)/i$. For any compact $K\subset(0,\infty)$ and integer $j\ge0$, the Leibniz rule is applied. The bounds from Lemma~\ref{lem:cutoff-derivatives} and Proposition~\ref{prop:fejer-bounds} (where $|\partial_x^{j-q} (F(x,i)/i)| \le D_{j-q} \cdot i$) yield a uniform majorant on $K$:
\[
\sup_{x\in K} |\partial_x^j h_i(x)| \le M_{j,K} \, i^{j+1} \, e^{-c_K i}.
\]
This is summable for any fixed $j$. By the Weierstrass M-test, all derivative series converge uniformly on compact sets, establishing that $\mathcal{P}_{\sigma}$ is a $C^\infty$ function.

For an integer $n$, Proposition~\ref{prop:fejer-properties} implies $F(n,i)/i = i \cdot \mathbf{1}_{i\mid n}$. The series reduces to a finite sum over divisors:
\[
\mathcal{P}_{\sigma}(n;\kappa)=\sum_{\substack{d\mid n\\ d\ge2}} d\,\phi_{\kappa}\!\Bigl(\frac{d}{n+1}\Bigr)-n.
\]
In the limit $\kappa\to\infty$, $\phi_{\kappa}(d/(n+1))\to 1$ for all divisors $d$. The sum becomes $\sum_{d\mid n, d\ge2} d = \sigma(n)-1$. The entire expression thus converges to $\sigma(n)-n-1$.
\end{proof}

\begin{conjecture}[Asymmetric companion zeros for $\mathcal{P}_{\sigma}$]\label{conj:Psigma-companions}
For each odd prime $p\ge 3$ and sufficiently large $\kappa$, it is conjectured that $\mathcal P_\sigma(\cdot;\kappa)$ displays an asymmetric pair of real zeros flanking $p$: an exponentially close left zero with
\[
p - x_p^-(\kappa)\ \asymp\ p\,\mathrm e^{-2\kappa/(p+1)}\qquad(\kappa\to\infty),
\]
and a right zero whose distance to $p$ appears to remain bounded away from zero as $\kappa\to\infty$. Numerics also indicate $\mathcal P_\sigma(p;\kappa) \approx -p\,\mathrm e^{-2\kappa/(p+1)}$.
\end{conjecture}

\begin{figure}[!htbp]
\centering
\includegraphics[width=0.9\textwidth]{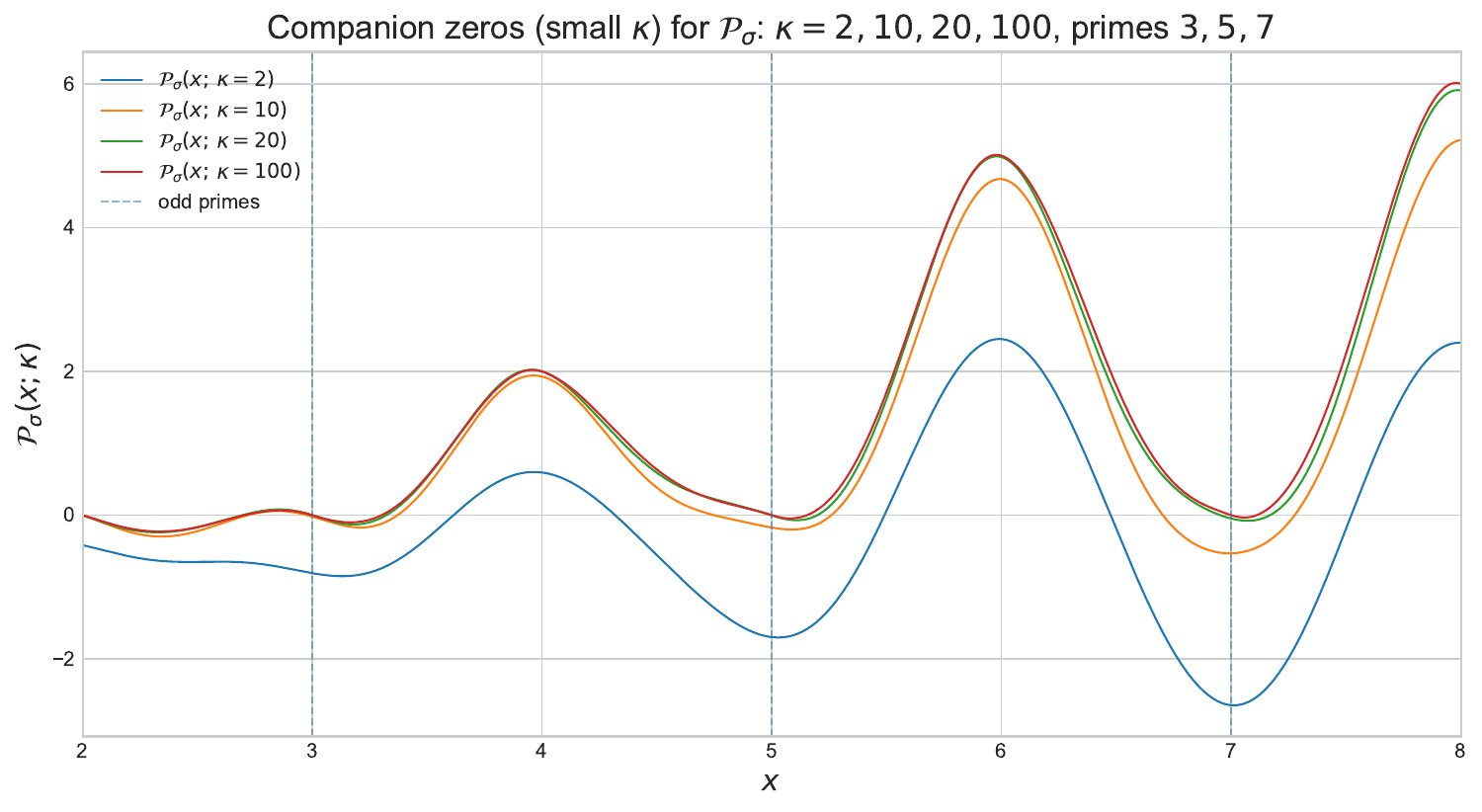}
\caption{Asymmetric companion zeros of $\mathcal{P}_{\sigma}(x;\kappa)$. Dashed lines mark primes $p=3,5,7$. Consistent with Conjecture~\ref{conj:Psigma-companions}, the left zero converges to $p$ exponentially fast with $\kappa$, while the right zero's distance to $p$ stabilizes at a positive limit. The value at $x=p$ is strictly negative for any finite $\kappa$.}
\label{fig:PsigmaCompanions}
\end{figure}

As a numerical illustration of Conjecture~\ref{conj:Psigma-companions}, 
Figure~\ref{fig:PsigmaCompanions} plots $\mathcal P_\sigma(x;\kappa)$ on $[2,8]$ for 
$\kappa\in\{2,10,20,100\}$. The predicted asymmetry is visible: an exponentially close left zero and a right zero at a $\kappa$-independent distance near each odd prime, while the value at $x=p$ remains negative for finite~$\kappa$.

\begin{lemma}[Uniform local derivative bounds near an odd prime]\label{lem:uniform-local-sigma}
Fix a prime $p\ge 3$. There exist $\varepsilon_0(p)>0$, $\kappa_0(p)\ge 1$, and $C_p>0$ such that, for all $\kappa\ge \kappa_0(p)$ and all $x$ with $|x-p|\le \varepsilon_0(p)$,
\[
\bigl|\partial_x^{m}\mathcal P_\sigma(x;\kappa)\bigr|\ \le\ C_p,\qquad m=2,3.
\]
In particular, the quadratic Taylor coefficient $A_p(\kappa)=\tfrac12\,\mathcal P_\sigma''(p;\kappa)$ is uniformly bounded below by $-C_p$, and the third-order remainder is bounded by $C_p|\varepsilon|^{3}$ uniformly in $\kappa\ge \kappa_0(p)$.
\end{lemma}

\begin{proof}
Write $\mathcal P_\sigma(x;\kappa)=\sum_{i\ge 2}\phi_\kappa(i/(x+1))\,F(x,i)/i - x$. Fix $p\ge 3$ prime and a compact interval $K=[p-\varepsilon_0,p+\varepsilon_0]\subset(0,\infty)$. By Proposition~\ref{prop:fejer-bounds}, $\partial_x^\ell(F(x,i)/i)\ll_\ell i$ uniformly in $x\in K$. By Lemma~\ref{lem:cutoff-derivatives},
\[
\sup_{x\in K}\bigl|\partial_x^q\phi_\kappa(i/(x+1))\bigr|\ \le\ C_{q,K,\kappa}\,i^q\,e^{-c_{K,\kappa} i},\qquad c_{K,\kappa}=\frac{2\kappa}{b+1},
\]
with $b=p+\varepsilon_0$. For $\kappa\ge \kappa_0(p)$ and any fixed $q$, the factor $C_{q,K,\kappa}$ grows at most polynomially in $\kappa$, whereas $e^{-c_{K,\kappa} i}$ decays exponentially in $i$ with rate proportional to $\kappa$. In particular,
\[
\sup_{\kappa\ge \kappa_0}\ C_{q,K,\kappa}\,i^q\,e^{-c_{K,\kappa} i}\ \le\ \widetilde C_{q,K,\kappa_0}\,i^q\,e^{-\widetilde c_K i},
\]
for some constants $\widetilde C_{q,K,\kappa_0}>0$ and $\widetilde c_K>0$ independent of $\kappa\ge\kappa_0$. By Leibniz' rule,
\[
\sup_{x\in K,\,\kappa\ge\kappa_0}\bigl|\partial_x^{m}\!\bigl(\phi_\kappa(i/(x+1))\,F(x,i)/i\bigr)\bigr|\ \ll_{m,K,\kappa_0}\ i^{m+1}\,e^{-\widetilde c_K i},
\]
which is summable in $i$ for $m=2,3$. Hence the derivative series converge uniformly on $K$, uniformly in $\kappa\ge\kappa_0$, by the Weierstrass $M$-test. Summation of the majorants yields the claimed uniform bounds for $\partial_x^{m}\mathcal P_\sigma(\cdot;\kappa)$ on $K$.
\end{proof}

\begin{lemma}[Asymptotic derivative at primes]\label{lem:Psigma-derivative-at-p}
Fix a prime $p\ge 3$. The derivative of the smooth divisor-sum analogue satisfies
\[
\lim_{\kappa\to\infty} \mathcal{P}_\sigma'(p;\kappa) = -1.
\]
\end{lemma}
\begin{proof}
Let $S_\sigma(x;\kappa) = \sum_{i\ge 2}\phi_\kappa(i/(x+1))\,F(x,i)/i$. Then $\mathcal{P}_\sigma'(x;\kappa) = S_\sigma'(x;\kappa)-1$. The derivative of the sum is evaluated termwise at $x=p$. The derivative of the $i$-th term is
\[
\partial_x\Bigl(\phi_\kappa\!\Bigl(\frac{i}{x+1}\Bigr)\frac{F(x,i)}{i}\Bigr) =
\Bigl(\partial_x\phi_\kappa\!\Bigl(\frac{i}{x+1}\Bigr)\Bigr)\frac{F(x,i)}{i} + \phi_\kappa\!\Bigl(\frac{i}{x+1}\Bigr)\frac{F'(x,i)}{i}.
\]
At $x=p$, one has $F(p,i)=0$ for any $i\nmid p$ and $F'(p,i)=0$ for any $i$ (cf.\ Proposition~\ref{prop:C1}). Since $p$ is prime, the only divisor in the sum is $i=p$. Thus, only the term for $i=p$ could be non-zero. For $i=p$, $F(p,p)=p^2$ and $F'(p,p)=0$. The sum reduces to a single term:
\[
S_\sigma'(p;\kappa) = \left. \partial_x\phi_\kappa\!\left(\frac{p}{x+1}\right)\right|_{x=p} \cdot \frac{F(p,p)}{p}
= p \cdot \left(-\frac{p}{(p+1)^2}\right) \phi_\kappa'\!\left(\frac{p}{p+1}\right).
\]
Since $\phi_\kappa'(u) = -\kappa\,\mathrm{sech}^2(\kappa(u-1))/2$, one has $\phi_\kappa'(p/(p+1)) \asymp -\kappa e^{-2\kappa/(p+1)}$, which vanishes as $\kappa\to\infty$. Therefore, $\lim_{\kappa\to\infty} S_\sigma'(p;\kappa) = 0$, and the result follows.
\end{proof}

\begin{lemma}[Existence and scale of the left zero for $\mathcal P_\sigma$]\label{lem:Psigma-left-zero}
Let $p\ge3$ be prime. There exist $\kappa_0(p)\ge1$ such that for all $\kappa\ge\kappa_0(p)$ the function $\mathcal P_\sigma(\cdot;\kappa)$ has a real zero $x_p^-(\kappa)$ to the left of $p$. With
\[
\Delta_p^{(\sigma)}(\kappa):=p\Bigl(1-\phi_\kappa\!\Bigl(\tfrac{p}{p+1}\Bigr)\Bigr)\asymp p\,\mathrm e^{-2\kappa/(p+1)},
\]
the distance to $p$ is bounded by
\[
\frac{1}{2}\,\Delta_p^{(\sigma)}(\kappa)\ \le\ p-x_p^-(\kappa)\ \le\ 2\,\Delta_p^{(\sigma)}(\kappa).
\]
\end{lemma}

\begin{proof}
Set $\Delta:=\Delta_p^{(\sigma)}(\kappa)$. By Lemma~\ref{lem:uniform-local-sigma}, a Taylor expansion at $x=p$ is valid:
\[
\mathcal P_\sigma(p-\varepsilon;\kappa)=\mathcal P_\sigma(p;\kappa) - \varepsilon\,\mathcal P_\sigma'(p;\kappa) + O(\varepsilon^2),
\]
where the $O(\cdot)$ term is uniform in $\kappa\ge\kappa_0(p)$. One has $\mathcal P_\sigma(p;\kappa)=-\Delta$. By Lemma~\ref{lem:Psigma-derivative-at-p}, $\mathcal P_\sigma'(p;\kappa)=-1+o(1)$ as $\kappa\to\infty$. The expansion becomes
\[
\mathcal P_\sigma(p-\varepsilon;\kappa)=-\Delta + \varepsilon\,(1+o(1)) + O(\varepsilon^2).
\]
For $\varepsilon_1=\Delta/2$, the value is $\mathcal P_\sigma(p-\varepsilon_1;\kappa)=-\Delta/2+o(\Delta)<0$ for sufficiently large $\kappa$.
For $\varepsilon_2=2\Delta$, the value is $\mathcal P_\sigma(p-\varepsilon_2;\kappa)=\Delta+o(\Delta)>0$ for sufficiently large $\kappa$.
By the Intermediate Value Theorem, a zero exists in the interval $(p-2\Delta, p-\Delta/2)$.
\end{proof}


\subsection{Critical discussion: non-integer zeros of $\mathcal{P}_{\sigma}$}
The indicator role of $\mathcal{P}_{\sigma}$ is confined to integer input. For non-integer $x$, the existence of zeros cannot be ruled out, because the function's sign is determined by the competition between the sum $S_{\sigma}(x;\kappa) = \sum_{i\ge2}\phi_{\kappa}(i/(x+1))\,F(x,i)/i$ and the linear term $x$. Near any integer, the factor $\sin^2(\pi x)$ in each $F(x,i)$ drives the sum to zero, ensuring $\mathcal{P}_{\sigma}(x;\kappa) \to -x$. Conversely, near rational points $x \approx im$ with small $i$, the term $F(x,i)/i$ can become large and lift $S_{\sigma}(x;\kappa)$ above the line $y=x$, creating intersections.

The following lemma provides a constructive lower bound for \emph{partial} sums of $S_{\sigma}$ away from such resonances. It does not yield a global sign for $\mathcal{P}_\sigma(x)=S_\sigma(x)-x$, as the two terms remain incomparable near integers and rational resonances. Consequently, no assertion about the sign of $\mathcal{P}_\sigma$ for non-integer arguments is claimed.

\begin{lemma}[Quantitative control away from resonances]\label{lem:constructive-lb}
Let $K=[a,b]\subset(0,\infty)$ be compact, $\eta\in(0,1)$, and $\kappa>0$. Choose $I_0\in\mathbb N$ and set $I:=\{2,\dots,I_0\}$. For $i\in I$ let
\[
E_i:=\{x\in K:\ \mathrm{dist}(x/i,\mathbb Z)<\delta_i\},\qquad 
\delta_i:=\frac{\eta}{8\,I_0\,\bigl((b-a)+I_0\bigr)}.
\]
Let $E_0:=\{x\in K:\ \mathrm{dist}(x,\mathbb Z)<\delta_0\}$ with $\delta_0:=\frac{\eta}{8(b-a+1)}$ and put $E:=E_0\cup\bigcup_{i\in I}E_i$. Then $\mathrm{meas}(E)\le \eta$ and, for all $x\in K\setminus E$,
\[
\sum_{i\in I}\phi_{\kappa}\!\Bigl(\frac{i}{x+1}\Bigr)\frac{F(x,i)}{i}
\ \ge\ 
\Bigl(\tfrac{2\delta_0}{\pi}\Bigr)^{\!2}\,
\sum_{i=2}^{I_0}\frac{\phi_{\kappa}\!\bigl(\tfrac{i}{a+1}\bigr)}{i}.
\]

\begin{proof}
For $x\in K\setminus E_0$, the distance $\mathrm{dist}(x,\mathbb Z)\ge \delta_0$ implies $\sin^2(\pi x)\ge (\tfrac{2\delta_0}{\pi})^2$. The map $x\mapsto \phi_\kappa(i/(x+1))$ is increasing on $[a,b]$ because $\phi_\kappa'(u)\le 0$ and $u'(x)<0$, hence
\[
\phi_\kappa\!\Bigl(\frac{i}{x+1}\Bigr)\ \ge\ \phi_\kappa\!\Bigl(\frac{i}{a+1}\Bigr)\qquad(x\in K).
\]
Using $F(x,i)=\sin^2(\pi x)/\sin^2(\pi x/i)$ and $\sin^2(\pi x/i)\le 1$ gives
\[
\phi_\kappa\!\Bigl(\frac{i}{x+1}\Bigr)\frac{F(x,i)}{i}
\ \ge\
\Bigl(\tfrac{2\delta_0}{\pi}\Bigr)^2\,\frac{\phi_\kappa\!\bigl(\tfrac{i}{a+1}\bigr)}{i}.
\]
Summing over $i\in I$ yields the lower bound. For the measure estimate, $\mathrm{meas}(E_0)\le 2\delta_0(b-a+1)$. For each $i\in I$, the condition $|x/i-m|<\delta_i$ is equivalent to $|x-im|<i\delta_i$, so $E_i$ is a union of at most $\lfloor(b-a)/i\rfloor+2$ intervals of length $2i\delta_i$. Therefore
\[
\mathrm{meas}\Bigl(\bigcup_{i=2}^{I_0}E_i\Bigr)\ \le\ \sum_{i=2}^{I_0}\Bigl(\Bigl\lfloor\frac{b-a}{i}\Bigr\rfloor+2\Bigr)\,2i\delta_i
\ \le\ 2\sum_{i=2}^{I_0}\bigl((b-a)+2i\bigr)\delta_i.
\]
With the chosen $\delta_i$, the sum is bounded by $\eta/2$, and together with $\mathrm{meas}(E_0)\le \eta/4$ this shows $\mathrm{meas}(E)\le \eta$.
\end{proof}
\end{lemma}

\begin{remark}[On the choice of bounds]
A lower bound for $\sin^2(\pi x)$ combined with the trivial upper bound $\sin^2(\pi x/i)\le 1$ yields a lower bound for $F(x,i)$. A lower bound on $\sin^2(\pi x/i)$ would only provide an upper bound for $F(x,i)$ and is therefore not useful here.
\end{remark}
\begin{remark}[On the choice of the thresholds $\delta_i$]
An $i$-dependent choice for $\delta_i$ could slightly improve the constants in the measure estimate. The uniform choice used above suffices for the proof and keeps the presentation simple.
\end{remark}

\begin{figure}[!htbp]
\centering
\includegraphics[width=0.9\textwidth]{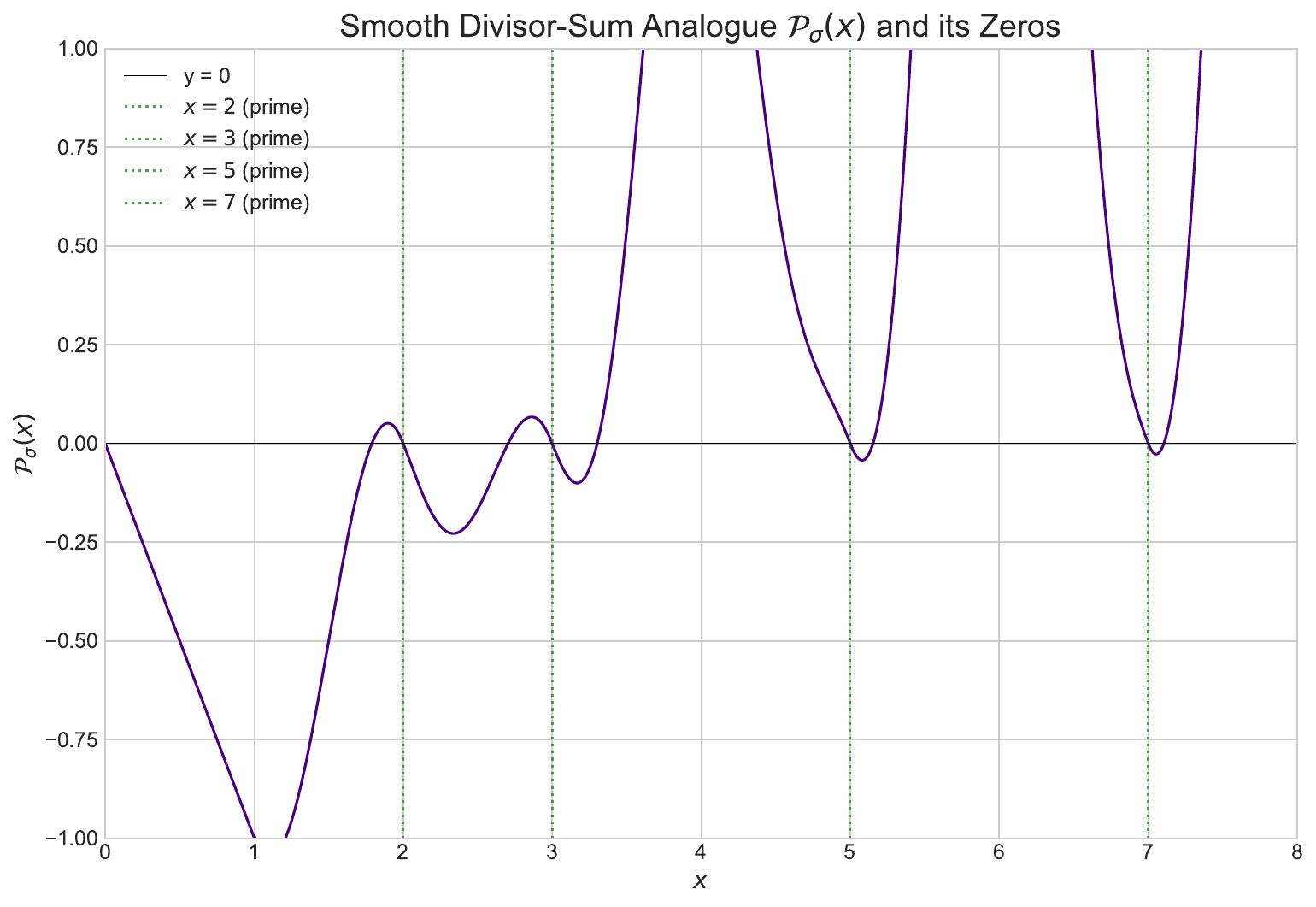}
\caption{Behavior of the smooth divisor-sum analogue $\mathcal{P}_{\sigma}(x; \kappa)$ for $\kappa=1000$ on $[0,8]$.
At any finite $\kappa$, $\mathcal P_\sigma(p;\kappa)<0$ at odd primes $p$; as $\kappa\to\infty$ the value at $x=p$ tends to $0$. Consistent with Conjecture~\ref{conj:Psigma-companions}, one observes a \emph{left prime-zero approximant} exponentially close to $p$ and a \emph{right companion zero} at a $\kappa$-independent distance. Away from these neighborhoods, the existence of additional non-integer zeros remains unresolved; accordingly, the prime-indicator interpretation is restricted to integer input.}
\label{fig:P_sigma}
\end{figure}

\FloatBarrier

\section{Concluding Remarks}

\subsection{Analytic Formulations and Arithmetic Consequences}
The argument uses the Fejér identity to connect trigonometric quotients with divisibility and evaluates at integers via the regularized cosine polynomial, while real inputs are handled by the quotient-of-sines representation. A consequential feature is that the Fejér-regularized quotient \(\bigl(\sin(\pi x)/\sin(\pi x/i)\bigr)^2\) yields, after normalization by \(i^2\), an exact divisor filter at integer arguments, \(\phi_i(n)=\mathbf{1}_{\,i\mid n}\). Superposition with weights then turns the filter into divisor sums—most notably \(\sum_i i^{2}\phi_i(n)=\sigma_{2}(n)\)—which provides the structural step enabling the smooth analogues \(\mathcal{P}_{\tau}\) and \(\mathcal{P}_{\sigma}\) via \(W(i)=i^{-2}\) and \(W(i)=i^{-1}\), with limits \(\mathcal{P}_{\tau}(n;\kappa)\to\tau(n)-2\) and \(\mathcal{P}_{\sigma}(n;\kappa)\to\sigma(n)-n-1\).

\subsection{Summary of Results and Open Questions}
Beyond a single indicator, smooth analogues are developed for arithmetic functions such as $\tau(x)$ and $\sigma(x)$. The function $\mathcal{P}_{\tau}$ is $C^\infty$ and approximates integer-prime zeros in the steep-cutoff limit. However, for finite steepness, it exhibits real zeros ("companion zeros") in the neighborhood of primes rather than at the primes themselves. The precise control and potential elimination of these analytic artefacts remains a key open question. For $\mathcal{P}_{\sigma}$, the indicator role is restricted to integer input, as the existence of non-integer zeros cannot be ruled out.

\subsection{Generalization to the Fej\'er--Dirichlet Lift}
The divisor-filter property, $F(n,i)/i^2 = \mathbf{1}_{i\mid n}$, is the central mechanism of this work. A subsequent manuscript~\cite{fuchs2025-fejer} generalizes this principle by constructing entire functions $\mathcal{T}_a(z) = \sum_{i\ge1} a(i) F(z,i)/i^2$ for arbitrary weight sequences $a(i)$. This operator, termed the \href{https://arxiv.org/abs/2509.12297}{Fej\'er--Dirichlet lift}, provides an analytic interpolant for the Dirichlet convolution $(a*1)(n)$. Its associated Dirichlet series is shown to factor as $\zeta(s)A(s)$, where $A(s)$ is the Dirichlet series of the weights. This framework reveals that the constructions herein are specific instances of a broader structural connection between Fejér analysis and the spectral theory of L-functions.

\subsection{Directions for Future Research}
Further directions include comparisons with other smooth prime indicators and alternative choices of $W(i)$ to model additional arithmetic functions, e.g., $\sigma_k(n)$ via $W(i)=i^{k-2}$. The properties of an analogous indicator for perfect numbers, such as $\mathcal{S}_\sigma(x;\kappa)-(2x-1)$, are left for future investigation. Recursive constructions based on indicator outputs (e.g., for $\omega(n)$) are another option; establishing their properties appears nontrivial. Practical efficiency for large inputs is not the focus here.

\appendix

\section{Illustrative prime-counting sums: constant-$C$ baseline and non-accumulative $H$-variant}

\noindent\textit{This appendix is illustrative. The constructions are not intended as practical prime-counting algorithms but serve to demonstrate how the integer prime-zero property of the smooth indicators can be utilized within summation formulas.}

\subsection{Construction Based on $\mathcal{P}_{\tau}$}
The construction is illustrative and not intended to compete with standard prime-counting methods.
The smooth analogue $\mathcal{P}_{\tau}(x)$ is used to construct related number-theoretic functions, rather than $\mathcal{P}(x)$.
The construction yields a prime–zero property for all integer inputs: for any integer $n \ge 2$, the limit $\lim_{\kappa\to\infty} \mathcal{P}_{\tau}(n;\kappa)=\tau(n)-2$ vanishes if and only if $n$ is prime.
This addresses the limitation that $\mathcal{P}(x)$ does not vanish at $p=2$.
This property allows a direct construction of a counting sum. The use of absolute values removes sign cancellations: for finite $\kappa$ the sign of $\mathcal{P}_{\tau}(x;\kappa)$ at non-integer $x$ is not controlled, whereas the intended feature of the summand is non-vanishing at non-primes.

\begin{definition}[Constant-threshold baseline]
An approximation for $\pi(x)$ for integers $x\ge2$ can be constructed via the sum
\begin{equation}\label{eq:pi_P_tau}
\pi_{\mathcal{P}_\tau}(x; C, \kappa) = \sum_{n=2}^{\lfloor x \rfloor} \left(1 - \frac{|\mathcal{P}_{\tau}(n; \kappa)|}{|\mathcal{P}_{\tau}(n; \kappa)| + C}\right),
\end{equation}
where $C$ is a small, positive constant (e.g., $C=0.001$) that serves as a fixed threshold.
The sum starts at $n=2$, the first prime, avoiding any ad-hoc adjustments.
\end{definition}

\begin{remark}[Choice of the threshold $C$]
A fixed $C$ controls the maximal per-term deviation for composites by $C/(1+C)$. For targets at the $10^{-3}$–level, $C=10^{-3}$ suffices. If an adaptive choice is preferred, setting $C=C(x):=\min\{10^{-3},\,(\log x)^{-2}\}$ preserves the qualitative behavior while reducing cumulative growth of the positive bias as $x$ increases.
\end{remark}

The behavior of the summand is as follows for large steepness parameter $\kappa$:
\begin{itemize}
    \item For prime $n$, $\mathcal{P}_{\tau}(n; \kappa)$ is close to $0$, so the summand is close to $1$.
    \item For composite $n$, $\mathcal{P}_{\tau}(n; \kappa)$ is close to $\tau(n)-2 \ge 1$, so the summand is close to $0$.
\end{itemize}
Consequently, for large $\kappa$ the summand in \eqref{eq:pi_P_tau} takes values close to $1$ at primes and close to $0$ at composites.

\subsection{Analysis of the Term-wise Error}
The deviation of this formula from the true value of $\pi(x)$ arises from two sources: the finite value of the steepness parameter $\kappa$, and the non-zero contribution of composite numbers controlled by the threshold $C$. For any finite $\kappa$, $\mathcal{P}_{\tau}(n; \kappa)$ deviates slightly from the integer limit $\tau(n)-2$.

Focusing on the systematic error from composite numbers in the ideal limit $\kappa\to\infty$, an upper bound for the error contribution of a single term can be established.

\begin{proposition}[Bounded Term-wise Error for the $\mathcal{P}_{\tau}$-based sum]
For any composite integer $n \ge 4$, the error contribution $E(n)$ is strictly positive and bounded by a constant dependent only on $C$. In the limit $\kappa\to\infty$, the error contribution is bounded above by
\[ E_{max} = \frac{C}{1+C}. \]
\end{proposition}

\begin{proof}
For a composite number $n$, the error contribution from the corresponding term in the sum is:
\[ E(n) = 1 - \frac{|\mathcal{P}_{\tau}(n; \kappa)|}{|\mathcal{P}_{\tau}(n; \kappa)| + C} = \frac{C}{|\mathcal{P}_{\tau}(n; \kappa)| + C}. \]
In the limit $\kappa\to\infty$, this becomes $E(n) = \frac{C}{|\tau(n)-2| + C}$.
The error is maximized when the denominator is minimized. This occurs when $|\tau(n)-2|$ is minimal. For any composite number $n$, the number of divisors $\tau(n)$ is at least 3 (e.g., for $n=p^2$). The smallest possible value for `n` is $4$, where $\tau(4)=3$. Thus, the minimum positive value of $|\tau(n)-2|$ for any composite $n$ is $|\tau(4)-2| = |3-2|=1$.
Therefore, the maximum possible error for any single composite term is bounded by:
\[ E_{max} = \frac{C}{1+C}. \]
For a choice of $C=0.001$, the maximum error for any single term is approximately $0.000999$.
\end{proof}

\subsection{Asymptotic Behavior of the Total Error}
Let $E(x;C,\kappa):=\pi_{\mathcal P_\tau}(x; C, \kappa) - \pi(x)$ and write the single-term deviation for composites as
\[
0< E(n)=1-\frac{|\mathcal P_\tau(n;\kappa)|}{|\mathcal P_\tau(n;\kappa)|+C}
=\frac{C}{|\mathcal P_\tau(n;\kappa)|+C}\qquad(n\ge4\ \text{composite}).
\]
(A) In the idealized limit $\kappa\to\infty$, one has $|\mathcal P_\tau(n;\infty)|=\tau(n)-2\ge1$, hence $E(n)\le \dfrac{C}{1+C}$ and
\[
0< E(x;C,\infty)\ \le\ \sum_{\substack{4\le n\le x\\ n\ \mathrm{composite}}}\frac{C}{1+C}
\ =\ \frac{C}{1+C}\,\bigl(\lfloor x\rfloor-\pi(x)-1\bigr)
\ =\ O(x).
\]
(B) For finite $\kappa$, the same reasoning yields bounds in terms of the uniform deviation
\[
\Delta_\kappa(x):=\max_{2\le n\le \lfloor x\rfloor}\Bigl|\mathcal P_\tau(n;\kappa)-\bigl(\tau(n)-2\bigr)\Bigr|,
\]
since $E(n)\le \dfrac{C}{\max\{1-\Delta_\kappa(x),0\}+C}$ for composites $n\le x$. Any effective control of $\Delta_\kappa(x)$ immediately transfers to a global bound for $E(x;C,\kappa)$.

\begin{proposition}[Global absolute deviation bounds for $\pi_{\mathcal P_\tau}$]
Let $C>0$ be fixed and $x\ge 2$.
\begin{itemize}
\item[(i)] For any $\kappa\in(0,\infty)$,
\[
\Bigl|\pi_{\mathcal P_\tau}(x;C,\kappa)-\pi(x)\Bigr|
\ \le\ \sum_{\substack{2\le n\le \lfloor x\rfloor\\ n\ \mathrm{composite}}} 1\;+\;\sum_{\substack{2\le n\le \lfloor x\rfloor\\ n\ \mathrm{prime}}} 1
\ =\ \lfloor x\rfloor-1
\ =\ O(x).
\]
Sharper $\kappa$-dependent bounds follow from the uniform deviation $\Delta_\kappa(x)$ introduced above.
\item[(ii)] In the idealized limit $\kappa\to\infty$,
\[
0\ \le\ \pi_{\mathcal P_\tau}(x;C,\infty)-\pi(x)
\ \le\ \frac{C}{1+C}\,\bigl(\lfloor x\rfloor-\pi(x)-1\bigr)
\ =\ O(x).
\]
\end{itemize}
\end{proposition}

\begin{proof}
Write $S(n)=1-\frac{|\mathcal P_\tau(n;\kappa)|}{|\mathcal P_\tau(n;\kappa)|+C}\in[0,1)$. Then
\[
\pi_{\mathcal P_\tau}(x;C,\kappa)-\pi(x)
=\sum_{\substack{n\le \lfloor x\rfloor\\ n\ \mathrm{comp}}}\!S(n)\;-\!\sum_{\substack{n\le \lfloor x\rfloor\\ n\ \mathrm{prime}}}\!\bigl(1-S(n)\bigr).
\]
By the triangle inequality and $S(n)\le 1$, one obtains (i). For (ii), as $\kappa\to\infty$, $S(n)\to 1$ for primes and $S(n)\to \frac{C}{\tau(n)-2+C}\le \frac{C}{1+C}$ for composites (since $\tau(n)-2\ge 1$). Hence the negative term vanishes and the positive term is bounded by the stated multiple of the number of composites up to $x$.
\end{proof}

This confirms the illustrative nature of the construction and explains why asymptotically sharp prime counting is not expected from \eqref{eq:pi_P_tau}. The formula is intended for illustrative and pedagogical purposes rather than for competition with classical prime-counting methods.

\subsection{A Variant with Summable Composite Leakage (H-Variant)}

To address the linear error accumulation from composite numbers in the baseline model, a variant is constructed where this leakage is bounded by a summable series. This property does not eliminate the error accumulation from primes (the "prime deficit"), which still grows with $\pi(x)$. Let $g(n):=|\mathcal P_\tau(n;\kappa(n))|$, where the steepness $\kappa(n):=\alpha(n+1)$ now depends on the input $n$ with a fixed rate $\alpha>0$. A dynamic threshold $\varepsilon(n):=(n+1)^{-\gamma}$ with $\gamma>1$ is introduced. The gating term is defined as
\[
a_n(\alpha,\gamma)\ :=\ \frac{\varepsilon(n)}{g(n)+\varepsilon(n)}\in(0,1],
\]
and the partial counting sum
\begin{equation}\label{eq:pi_H}
\pi_H(x;\alpha,\gamma)\ :=\ \sum_{n=2}^{\lfloor x\rfloor} a_n(\alpha,\gamma).
\end{equation}

\noindent\emph{Terminology.} The term “non-accumulative” refers exclusively to the composite leakage
\(
\sum_{n\ \mathrm{composite}} \varepsilon(n)/(g(n)+\varepsilon(n)),
\)
which is summable for $\gamma>1$. The prime deficit
\(
\sum_{p\le x} B(\alpha)/(B(\alpha)+\varepsilon(p))
\)
accumulates with $\pi(x)$.

For any odd prime $p$, the argument of the cutoff function $\phi_{\kappa(p)}$ becomes independent of $p$. With $\kappa(p)=\alpha(p+1)$, the term $u-1$ inside the hyperbolic tangent evaluates to $\frac{p}{p+1}-1 = -\frac{1}{p+1}$, which yields $\kappa(p)(u-1)=-\alpha$. Consequently, $g(p)$ simplifies to a constant that depends only on $\alpha$:
\[
g(p) = \left|\phi_{\kappa(p)}\!\left(\frac{p}{p+1}\right)-1\right| = \left|\frac{1-\tanh(-\alpha)}{2}-1\right| = \frac{1-\tanh(\alpha)}{2} = \frac{1}{e^{2\alpha}+1}.
\]
This constant is denoted by $B(\alpha)$. The gating term $a_p$ is thus not approximately 1, but rather a specific value determined by the competition between the prime residual $B(\alpha)$ and the threshold $\varepsilon(p)$. The mechanism is effective only if $\varepsilon(p) \gg B(\alpha)$, ensuring $a_p$ is close to 1. For composite $n$, $g(n)$ is expected to be of order 1.

\paragraph*{Finite-range admissibility for the non-accumulative $H$-variant}
Let $\kappa(n)=\alpha(n+1)$ with $\alpha>0$, $\varepsilon(n)=(n+1)^{-\gamma}$ with $\gamma>1$, and
\[
t_n=\frac{\varepsilon(n)}{\,|\mathcal P_\tau(n;\kappa(n))|+\varepsilon(n)}\,.
\]
Define $B(\alpha):=\dfrac{1-\tanh\alpha}{2}=\dfrac{1}{e^{2\alpha}+1}$.

\begin{lemma}[Prime residual at odd primes]
For every odd prime $p$,
\[
\bigl|\mathcal P_\tau(p;\kappa(p))\bigr| \;=\; B(\alpha)=\frac{1}{e^{2\alpha}+1}.
\]
\end{lemma}

\begin{lemma}[Uniform composite gap]\label{lem:uniform-composite-gap}
Let $\kappa(n)=\alpha(n+1)$ with $\alpha>0$. For every composite integer $n\ge 4$,
\[
\mathcal P_\tau(n;\kappa(n))\;\ge\; 1-2\,e^{-2\alpha}\;=:\;c_\alpha.
\]
In particular $c_\alpha>0$ holds as soon as $\alpha>\tfrac12\log 2$.

\begin{proof}
For $n\ge 2$,
\[
\mathcal P_\tau(n;\kappa(n))=\sum_{\substack{d\mid n\\ d\ge 2}}\phi_{\kappa(n)}\!\Bigl(\frac{d}{n+1}\Bigr)-1.
\]
For any divisor $d\le n$,
\[
\kappa(n)\Bigl(\frac{d}{n+1}-1\Bigr)=\alpha\,(d-n-1)\le -\alpha,
\]
hence
\[
\phi_{\kappa(n)}\!\Bigl(\tfrac{d}{n+1}\Bigr)
=\frac{1-\tanh\!\bigl(\alpha(d-n-1)\bigr)}{2}
=\frac{1+\tanh\!\bigl(\alpha(n+1-d)\bigr)}{2}
\ \ge\ 1-e^{-2\alpha},
\]
since $\tanh t\ge 1-2e^{-2t}$ for $t\ge 0$. Because $n$ is composite, there exists a proper divisor $d$ with $2\le d\le n/2$. For this $d$ and for $d=n$ one has $n+1-d\ge 1$, hence the displayed lower bound applies to both indices. Therefore the sum over $d\ge 2$ contains at least two terms not smaller than $1-e^{-2\alpha}$, and
\[
\mathcal P_\tau(n;\kappa(n))
\ \ge\ 2(1-e^{-2\alpha})-1
\ =\ 1-2e^{-2\alpha}.
\]
\end{proof}
\end{lemma}

\begin{corollary}[Finite-range admissibility]
If $\varepsilon(X)\ge \lambda\,B(\alpha)$ for some $\lambda\ge 1$, then $t_p\ge \dfrac{\lambda}{1+\lambda}$ for all primes $p\le X$.
For $\varepsilon(n)=(n+1)^{-\gamma}$ this is ensured by
\[
\gamma \;\le\; \frac{\log\!\bigl(1/B(\alpha)\bigr)-\log \lambda}{\log(X+1)}
 \;=\; \frac{\log\!\bigl(e^{2\alpha}+1\bigr)-\log \lambda}{\log(X+1)}.
\]
\end{corollary}

\paragraph*{Parameter choice for the figure}
With $\alpha=18.5$, the prime residual is $B(18.5)=(e^{37}+1)^{-1}$. For the range $X=50$, the condition $\varepsilon(X) \ge \lambda B(\alpha)$ with $\lambda=100$ requires $(51)^{-\gamma} \ge 100 \cdot (e^{37}+1)^{-1}$. This yields an upper bound on $\gamma$:
\[
\gamma \le \frac{\log(e^{37}+1) - \log(100)}{\log(51)} \approx \frac{37.00 - 4.61}{3.93} \approx 8.24.
\]
Here $\log$ denotes the natural logarithm. The choice $\gamma=5$ is admissible under this constraint and ensures $t_p \ge \lambda/(1+\lambda) = 100/101$ for all primes $p\le 50$.

\begin{proposition}[Error decomposition for the H-variant]
The absolute error satisfies
\[
\bigl| \pi_H(x;\alpha,\gamma) - \pi(x) \bigr|
\;\le\; E_c(x) + E_p(x),
\qquad
E_c(x):=\sum_{\substack{n\le \lfloor x\rfloor \\ n\ \mathrm{composite}}}\frac{\varepsilon(n)}{g(n)+\varepsilon(n)},
\quad
E_p(x):=\sum_{p\le x}\frac{B(\alpha)}{B(\alpha)+\varepsilon(p)}.
\]
With $g(n):=\bigl|\mathcal P_\tau(n;\kappa(n))\bigr|$ and Lemma~\ref{lem:uniform-composite-gap},
\[
g(n)\;\ge\; c_\alpha:=1-2e^{-2\alpha}\qquad(n\ \mathrm{composite}),
\]
and therefore, for every $\gamma>1$,
\[
E_c(x)\ \le\ \frac{1}{c_\alpha}\sum_{n=4}^{\infty}(n+1)^{-\gamma}\ <\ \frac{\zeta(\gamma)}{c_\alpha}\,.
\]
In particular, $E_c(x)$ is uniformly bounded in $x$ provided $\alpha>\tfrac12\log 2$.
\end{proposition}

\begin{proof}[Proof sketch]
The total error is the sum over composites (leakage) minus the sum of deviations from 1 at primes (deficit). The absolute error is bounded by the sum of their magnitudes via the triangle inequality.

For the composite leakage $E_c(x)$, each term $a_n$ is bounded above by $\varepsilon(n)/c_{\alpha,x}$. Since the series $\sum \varepsilon(n)$ converges for $\gamma>1$, the total leakage is bounded by a constant that does not grow with $x$.

For the prime deficit $E_p(x)$, the deviation at each prime is $1-a_p = B(\alpha)/(B(\alpha)+\varepsilon(p))$. This term does not form a summable series over the primes. While each term is small for a suitable choice of parameters, their sum accumulates with $\pi(x)$. The "non-accumulative" property of the H-variant thus applies only to the error contribution from composite numbers.
\end{proof}

\begin{remark}[Quantitative separation on practical ranges]
For every divisor $d\mid n$ with $d\le n$,
\[
\phi_{\kappa(n)}\!\Bigl(\frac{d}{n+1}\Bigr)\ \ge\ 1-e^{-2\alpha},
\]
hence Lemma~\ref{lem:uniform-composite-gap} gives the uniform bound
$g(n)\ge c_\alpha=1-2e^{-2\alpha}$ for composites $n$. Consequently,
$E_c(x)\le c_\alpha^{-1}\sum_{n\ge 4}(n+1)^{-\gamma}$ is non-accumulative for any $\gamma>1$ and every $\alpha>\tfrac12\log 2$.
Empirically, with $(\alpha,\gamma)=(19,7)$ on $x\le 10^4$ the observed deviation was $\approx 1.4\times 10^{-5}$ in double precision, while the constant-$C$ baseline shows a linear drift.
\end{remark}

\begin{figure}[!htbp]
\centering
\includegraphics[width=0.9\textwidth]{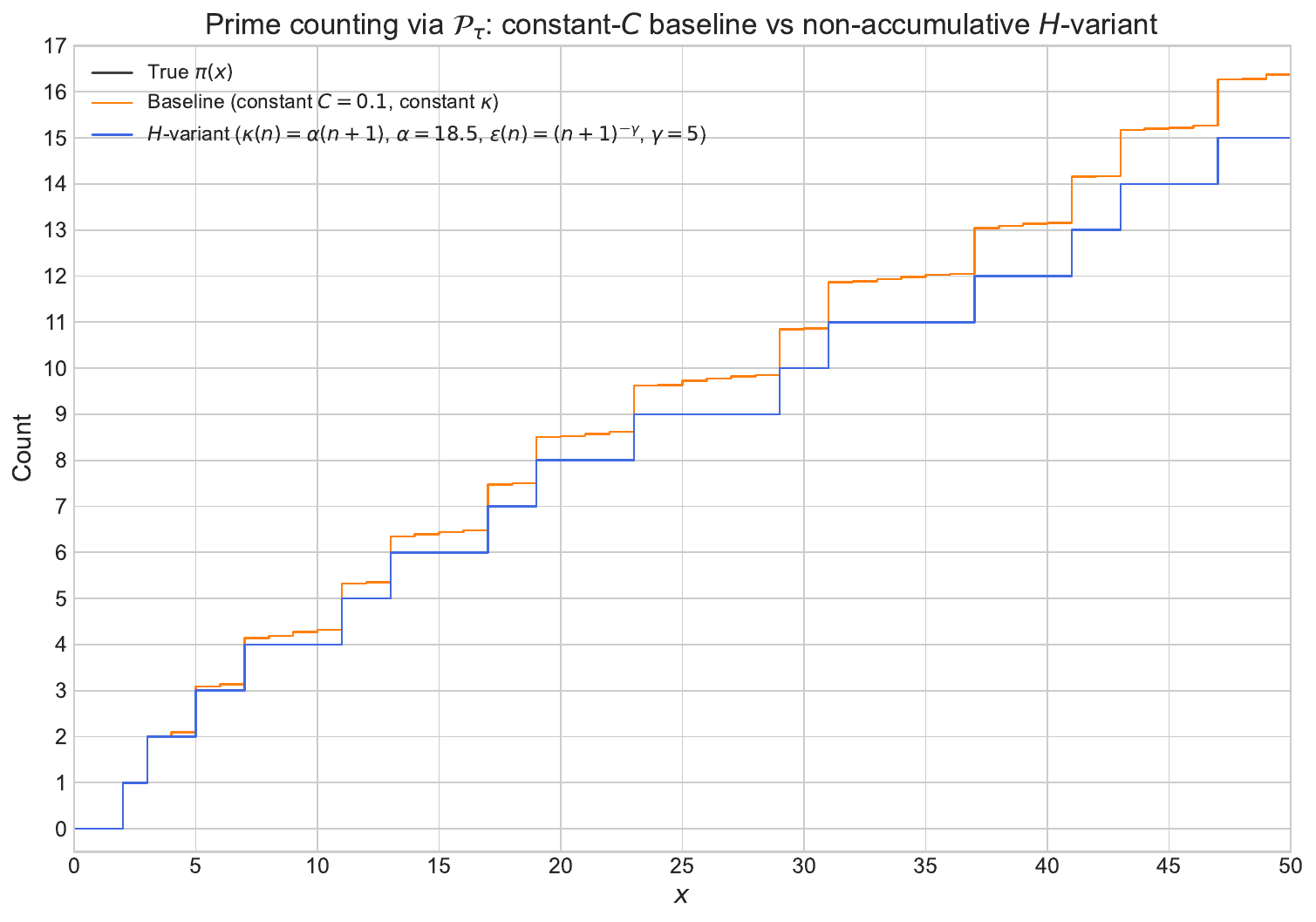}
\caption{Comparison of two illustrative prime-counting sums on $0\le x\le 50$. The constant-threshold baseline uses $C=0.1$ and a fixed steepness $\kappa$, which produces a visible linear drift due to cumulative composite leakage. The non-accumulative $H$-variant uses $\kappa(n)=\alpha(n+1)$ and $\varepsilon(n)=(n+1)^{-\gamma}$ with $\alpha=18.5$ and $\gamma=5$, yielding a curve that visually coincides with the staircase of $\pi(x)$ on this range. Both constructions are illustrative; no asymptotic claims are made.}
\label{fig:primecounting}
\end{figure}
\FloatBarrier

\section{Reference implementation}\label{app:code}

\noindent\textit{Algorithm for stable evaluation of $\mathcal{P}(x)$.}
\begin{lstlisting}[language=Python]
# Numerically stable O(sqrt(x)) evaluation of P(x) with explicit resonance guards (hardened).

from math import sin, cos, pi, sqrt, ceil

def P(x, eps=1e-12):
    assert x > 0.0
    N = ceil(sqrt(x))
    S = 0.0

    # Global tolerances
    EPS_INT = eps
    EPS_RES = 1e-6

    # Reuse numerator robustly
    sx = sin(pi * x)

    # Integer detection
    n = int(round(x))
    near_int = abs(x - n) < EPS_INT

    def is_resonant(t, i):
        # t ~ integer or sin(pi t) small -> potential denominator cancellation
        return (abs(t - round(t)) < max(EPS_RES, 0.1 / i)) or (abs(sin(pi * t)) < EPS_RES)

    for i in range(2, N + 1):
        t = x / i

        if is_resonant(t, i):
            # --- Resonance detected: x/i is close to an integer ---
            m = int(round(t))
            delta = x - i * m  # Deviation from the nearest integer multiple of i

            # Exact integer and exact divisibility handled purely in integer arithmetic
            if near_int and (n % i == 0):
                S += i * i    # Exact value F(n,i) = i^2 for divisors
                continue

            # Local even Taylor surrogate near resonance
            if abs(delta) <= min(i / 8.0, 1.0 / pi):
                alpha = (1.0 / 6.0) * (1.0 - 1.0 / (i * i))
                val = i * i * (1.0 - 2.0 * alpha * (pi * delta) ** 2)
                S += val if val > 0.0 else 0.0
            else:
                # Fall back to cosine polynomial to avoid cancellation
                A = 0.0
                for k in range(1, i):
                    A += (i - k) * cos(2 * pi * k * x / i)
                S += i + 2.0 * A
        else:
            denom = sin(pi * t)
            S += (sx * sx) / (denom * denom)

    return S / x
\end{lstlisting}

\paragraph*{Implementation notes: resonance guards and control flow}
Note that Listing 1 serves as a conceptual illustration of the complete stabilization logic; for modularity and clarity, the companion script implements this logic across the main function and a dedicated \texttt{fejer\_closed\_form\_term} helper. The routine separates indices into non-resonant and resonant regimes and applies the least sensitive representation in each case.

\begin{itemize}
\item \emph{Resonance detection.} For $t=x/i$ and $m=\mathrm{round}(t)$, a resonance is flagged if $\lvert t-m\rvert$ is small or $\lvert\sin(\pi t)\rvert$ is small. This captures proximity to the pole set $x\approx im$ where the denominator of $\sin(\pi x/i)$ is tiny.

\item \emph{Exact divisors.} If $x$ is (numerically) an integer and $i\mid x$, the value $F(x,i)=i^2$ is returned, which avoids cancellation.

\item \emph{Local $\delta$-switch near resonance.} Writing $x=im+\delta$, the local even expansion from Lemma~\ref{lem:resonance} yields
\[
\Bigl(\tfrac{\sin(\pi x)}{\sin(\pi x/i)}\Bigr)^{2}
= i^{2}\bigl(1-2\alpha_i(\pi\delta)^{2}\bigr)+O\!\bigl(i^{2}(\pi\delta)^{4}\bigr),
\qquad \alpha_i=\tfrac16\bigl(1-\tfrac{1}{i^{2}}\bigr).
\]
The code uses this quadratic truncation for $\lvert\delta\rvert$ below a fixed threshold (chosen to keep the remainder within the fourth-order scale guaranteed by Lemma~\ref{lem:resonance}) and clamps the truncated value at $0$ only within this Taylor branch. This ensures $F\ge 0$ and induces a one-sided (downward) truncation error confined to the local surrogate.

\item \emph{Safe fallback away from the local window.} If $\lvert\delta\rvert$ exceeds the local window, the cosine-polynomial form
\[
F(x,i)=i+2\sum_{k=1}^{i-1}(i-k)\cos\!\Bigl(\frac{2\pi k x}{i}\Bigr)
\]
is used. This avoids loss of significance from nearly cancelling sines.

\item \emph{Non-resonant path.} When no resonance is detected, the quotient-of-sines formula
\[
F(x,i)=\frac{\sin^{2}(\pi x)}{\sin^{2}(\pi x/i)}
\]
is numerically stable and is evaluated directly, reusing $\sin(\pi x)$.

\item \emph{Cost and coverage.} The non-resonant path costs $O(1)$ per index and dominates $O(\sqrt{x})$ indices. The local surrogate and the polynomial fallback are triggered only on a small fraction of indices near rational resonances; by Lemma~\ref{lem:resonance}, the total measure of such indices is controlled, so the overall cost remains $O(\sqrt{x})$.

\item \emph{Parameter choices.} The thresholds $\lvert t-m\rvert$, $\lvert\sin(\pi t)\rvert$, and the admissible $\lvert\delta\rvert$ window are set to keep the fourth-order remainder within the scale specified in Lemma~\ref{lem:resonance}. These can be tightened for higher accuracy without changing the control flow.
\end{itemize}

\smallskip
Approach (B) has $O(\sqrt{x})$ complexity and is preferable when run time dominates, while approach (A) is suited to verification and the derivation of analytic identities.

\paragraph*{Acknowledgements}
Valuable discussions and helpful feedback from colleagues are gratefully acknowledged.

\vfill
\par\noindent
\small
Sebastian Fuchs \\[1ex]
\textit{arXiv:} \texttt{\href{https://arxiv.org/abs/2506.18933}{arXiv:2506.18933}} \\
\textit{DOI:} \texttt{\href{https://doi.org/10.5281/zenodo.17360383}{10.5281/zenodo.17360383}} \\ 
\textit{ORCID:} \texttt{\href{https://orcid.org/0009-0009-1237-4804}{0009-0009-1237-4804}} \\

\end{document}